\documentclass[final,3p,times,12pt]{elsarticle}
\usepackage[fleqn]{amsmath}
\usepackage{amsthm}
\usepackage{mathrsfs}
\usepackage{amssymb}
\usepackage{graphicx}
\usepackage{subfigure}
\usepackage{color}
\usepackage{array}
\usepackage{mathrsfs}
\usepackage{float}
\usepackage{epstopdf}
\usepackage{booktabs}
\usepackage{tabularx}
\usepackage{lineno}
\usepackage{grffile}
\usepackage{extarrows}
\usepackage{microtype}
\usepackage{ragged2e}
\allowdisplaybreaks[4]
\graphicspath{{Fig/}}
\usepackage[colorlinks=true, linkcolor=blue, citecolor=blue, urlcolor=blue]{hyperref}
\biboptions{sort&compress}

\newcommand{\E}{\mathbb{E}}
\newcommand{\PP}{\mathbb{P}}
\newcommand{\RR}{\mathbb{R}}

\def\[{\langle} \def\]{\rangle}
\newdefinition{lemma}{Lemma}[section]
\newdefinition{theorem}{Theorem}[section]
\newdefinition{definition}{Definition}[section]
\newdefinition{assumption}{Assumption}[section]
\newdefinition{proposition}{Proposition}[section]
\newdefinition{remark}{Remark}[section]
\newdefinition{claim}{Claim}[section]
\newdefinition{corollary}{Corollary}[section]
\newdefinition{example}{Example}[section]

\begin{document}
	\newpage
	\begin{frontmatter}
		\setcounter{page}{1}
		\title{An explicit domain-preserving projection scheme for scalar stochastic differential equations}
		\author{Xingwei Hu$^{a,b}$,~~Yaozhong Hu$^{b}$,~~Aiguo Xiao$^{a,c,*}$}
			\footnote{\textsuperscript{$*$} Corresponding author (Aiguo Xiao).\\
			\emph{Email addresses}: \texttt{xingweihu@smail.xtu.edu.cn} (X.\ Hu), 
			\texttt{yaozhong@ualberta.ca} (Y.\ Hu),
			\texttt{xag@xtu.edu.cn} (A.\ Xiao).}
		\address{$^a$School of Mathematics and Computational Science $\&$ Hunan Key Laboratory for Computation and Simulation in Science and Engineering, Xiangtan University, Xiangtan, Hunan 411105, China \\
			$^b$ Department of Mathematical and Statistical Sciences, University of Alberta, Edmonton T6G 2G1, Alberta, Canada\\
            $^c$ Hunan Research Center of the Basic Discipline Fundamental Algorithmic Theory and Novel Computational Methods $\&$ China National Center for Applied Mathematics in Hunan, Xiangtan University, Xiangtan, Hunan 411105, China }
		\date{}

		\begin{abstract}
			\par We develop an explicit domain-preserving projection scheme for scalar stochastic differential equations whose solutions evolve in a prescribed interval. A Lamperti-type transformation converts the original equation into one with additive noise, whose drift may be singular near finite boundaries or grow polynomially at infinity. The proposed domain-preserving projected Euler--Maruyama (DPPEM) scheme combines an explicit Euler--Maruyama step with projection onto a step-size-dependent interior interval. Under one-sided Lipschitz, boundary, and weighted growth assumptions, with projection parameters chosen according to the coefficient growth and available moment orders, we establish strong order one convergence for admissible moment orders. The analysis combines local consistency and stability estimates with moment and boundary negative moment bounds for the exact solution. Numerical moment bounds follow from the convergence estimate, while numerical negative moment bounds are obtained under additional conditions. Estimates for the inverse Lamperti map yield strong order one convergence for the original equations in the CIR, Heston 3/2, CEV, A\"it–Sahalia, and bounded interest-rate models under the stated model conditions. Numerical experiments illustrate the predicted convergence rates and compare computational costs with the Lamperti backward Euler–Maruyama scheme.
		\end{abstract}
		
	\end{frontmatter}
	
\section{Introduction}
	\label{sec.1}
    
	Stochastic differential equations (SDEs) evolving in a prescribed interval arise naturally in mathematical finance, population dynamics, epidemiology, and many other applications. Typical domains include the positive half-line and bounded intervals, with positivity or boundedness forming an intrinsic part of the SDE. The coefficients may fail to be globally Lipschitz, degenerate near the boundary, or exhibit superlinear growth at infinity. After a Lamperti-type transformation, such boundary degeneracy may further lead to singular drift behavior near the transformed boundaries. These features complicate both the construction of numerical schemes and their convergence analysis.
    
    Standard explicit schemes generally do not preserve the transformed domain, as an Euler--Maruyama (EM) or Milstein step may leave the prescribed interval with positive probability. This limitation has motivated extensive research on domain-preserving numerical approximations for scalar SDEs. For the stochastic Susceptible-Infected-Susceptible model \cite{GGHM2011}, several explicit schemes with strong order one have been developed \cite{YH2021, CGW2021, LW2023}. Positivity-preserving discretizations have also been developed for the constant elasticity of variance (CEV) model \cite{KAHL2008,HS2015}, the Heston \(3/2\) model \cite{WG2025}, and the A\"it--Sahalia (AS) model \cite{DFWM2023,LCW2025,JLW2025}. The Cox--Ingersoll--Ross (CIR) process \cite{CIR1985} has been extensively investigated by both implicit and explicit approaches \cite{AA2005,SNS2012,AA2013,AC2020}. 
    
    General approaches to scalar SDEs on intervals use Lamperti or logarithmic transformations, implicit and explicit schemes, truncation, and projection \cite{LNX2025,MS2007,NS2014,brc26,EL2026,PPlogTM,TX2024,LG2023,HXW}. Related positivity-preserving schemes for multidimensional models and stochastic systems are discussed in \cite{Niko2015,zhou_hu_zhao,hxw2026,LC2023,CWF2024,HC2025,THX2026}. For scalar equations with positive diffusion coefficients, a Lamperti-type transformation yields additive noise but may introduce drift singularities at finite endpoints. Under suitable assumptions, implicit schemes can yield well-defined, domain-preserving approximations for transformed equations with singular drifts near the boundary \cite{NS2014,LNX2025}. For explicit methods using step-size-dependent truncation or projection, convergence analysis must control the errors introduced by these modifications.

   Beyn, Isaak, and Kruse \cite{BIK2016,BIK2017} developed a framework for strong convergence analysis based on stochastic C-stability and B-consistency and applied it to projected EM and Milstein schemes. Their analysis accommodates superlinear coefficient growth under global monotonicity and polynomial regularity assumptions, without requiring a priori higher moment bounds for the numerical approximation. Their schemes first project the current approximation onto a ball whose radius increases as the step size decreases and then perform the update using the projected state. For additive noise, their results yield strong order one convergence of the projected EM scheme \cite{BIK2017}. These results assume coefficients defined on $\mathbb R^d$ and do not directly cover drifts singular at finite interval endpoints. Interval problems thus require additional control near finite boundaries and, when the domain is unbounded, at infinity.

   To account for drift singularities near finite endpoints and polynomial growth at infinity, we propose a fully explicit DPPEM scheme and provide a unified analysis of its strong convergence. Consider the scalar SDE
    \begin{equation}\label{sde}
    \mathrm{d}U_t = F(U_t)\,\mathrm{d}t + G(U_t)\,\mathrm{d}B_t,
    \quad U_0=u_0\in(A,B),
    \end{equation}
    where $-\infty\le A < B\le\infty$ and $G(x)>0$ for $x\in(A,B)$. We apply a Lamperti-type transformation $W_t=\Phi(U_t)$, which converts \eqref{sde} into an SDE with additive noise,
    \begin{equation*}
\mathrm{d}W_t=\mu(W_t)\,\mathrm{d}t+\sigma\,\mathrm{d}B_t,\quad W_0\in(a,b),
    \end{equation*}
    where $\sigma>0$ and \((a,b)\) may have finite or infinite endpoints. The coefficient assumptions concern the transformed drift $\mu$. To ensure well-defined drift evaluations and control coefficient growth, we introduce a step-size-dependent closed interior interval $\mathbb D_h=[a_h,b_h]\subset(a,b)$. For sufficiently small $h>0$, $\mathbb D_h$ contains $W_0$, and its endpoints approach those of $(a,b)$ as $h\rightarrow0$. The DPPEM scheme is
\begin{equation*}
W^{n+1}=\Pi_{\mathbb D_h}\bigl(W^n+h\mu(W^n)+\sigma\Delta B_n\bigr),
\quad W^0=W_0,
\end{equation*}
where $\Pi_{\mathbb D_h}$ is the Euclidean projection onto $\mathbb D_h$. Since interval projection is scalar clipping, the scheme is fully explicit and satisfies $W^n\in\mathbb D_h$ at every time step, pathwise. This keeps numerical values a positive distance from finite endpoints and imposes finite bounds on unbounded domains.

     The exact solution may leave $\mathbb D_h$ while remaining in $(a,b)$. Enlarging $\mathbb D_h$ reduces this probability but may increase the drift and derivative bounds within $\mathbb D_h$. To achieve strong order-one convergence, the projection parameters must therefore balance the coefficient bounds needed for consistency and stability with the error contributions on exit events. Under the regularity, one-sided Lipschitz, and boundary assumptions, we first establish global well-posedness and moment and boundary negative moment estimates for the transformed equation. Moment estimates control excursions beyond the finite cutoffs on unbounded intervals, while boundary negative moment estimates control excursions near finite endpoints. 
     
     The weighted growth conditions in Assumption~3.1 distinguish the requirements for local consistency from those for stability. Combining these conditions with the exact-solution estimates gives explicit restrictions linking the projection parameters to coefficient growth and the available moment and negative moment orders. We combine a stopped error estimate with bounds on the probability that the exact solution leaves $\mathbb D_h$ at a grid point to obtain global error control. The argument treats finite and infinite endpoints without assuming numerical moment or negative moment bounds uniform in $h$ a priori.

     Under these restrictions, we establish strong order-one convergence in the $L^r$ norm of the maximum error over the time grid for admissible $r>1$. Numerical moment bounds follow for admissible orders, while numerical boundary negative moment bounds require additional conditions. For representative models, estimates for the inverse Lamperti map transfer convergence to the original variables under model specific parameter restrictions. Numerical experiments illustrate the convergence rates and compare computational costs.
    
    The main contributions are as follows.
\begin{itemize}
	\item We establish global well-posedness together with moment and negative moment  estimates for the transformed SDE under the regularity, one-sided Lipschitz, and boundary assumptions.
	\item We construct a fully explicit DPPEM scheme and provide a unified strong convergence analysis allowing drift singularities near finite endpoints and polynomial growth at infinity. Under the stated coefficient assumptions and projection parameter restrictions, we establish strong order one convergence for the transformed equation in the $L^r$ norm of the maximum error over the time grid for admissible $r>1$.
	\item We specify admissible projection parameters by combining separate growth bounds for consistency and stability with exact-solution moment and boundary negative moment estimates. The analysis covers finite and infinite endpoints without assuming numerical moment or boundary negative moment bounds uniform in $h$ a priori. Numerical moment bounds follow for admissible orders, while numerical boundary negative moment bounds require additional conditions.
\end{itemize}
    
    The remainder of the paper is organized as follows. Section 2 introduces notation and preliminaries. Section~3 develops the DPPEM scheme and establishes its strong order one convergence. Section 4 applies the theoretical results to representative scalar models. Section 5 reports numerical experiments, and Section 6 concludes the paper.

	\section{Notation and preliminaries }
	\label{sec.2}
   Throughout this paper, let
$(\Omega,\mathcal F,\{\mathcal F_t\}_{t\in[0,T]},\mathbb P)$
be a filtered probability space satisfying the usual conditions. We denote expectation by $\mathbb E$. For an event or a set $A$, $\mathbb I_A$ denotes its
indicator, and $A^c$ denotes its complement.
For a random variable $X$ and $p\ge1$, we write $\|X\|_{L^p}:=(\mathbb E|X|^p)^{1/p}$. For a set $A$, $\overline A$ denotes its closure, and
$A\Subset D$ means that $\overline A$ is a compact subset of $D$.
If $K\subset\mathbb R^d$ is a nonempty closed convex set,
$\Pi_K$ denotes the Euclidean projection onto $K$, i.e., $\Pi_K(x):=\arg\min_{y\in K}|x-y|$. In particular, $\Pi_K$ is nonexpansive: $|\Pi_K(x)-\Pi_K(y)|\le |x-y|$. For real numbers $m,n$, we write $m\wedge n:=\min\{m,n\}$ and
$m\vee n:=\max\{m,n\}$. Let $N\in\mathbb N$, $h=T/N$, and $t_n:=nh$, $n=0,\ldots,N$.
We set $\mathbb E_n[\cdot]
:=\mathbb E[\cdot\mid\mathcal F_{t_n}]$. For $\varphi \in C^2((a,b);\mathbb R)$, define the generator
	\begin{equation*}
		\mathcal{L}\varphi(x) := \mu(x)\varphi'(x) + \frac{\sigma^2}{2}\varphi''(x).
	\end{equation*}
    Throughout the paper, $C$ denotes a generic positive constant
independent of $h$ and $n$, whose value may change from line to line. 
	
	\begin{proposition}[Ville's inequality \cite{JOb2015,HR2020}]\label{ville}
		Suppose $Z_t$ is a nonnegative supermartingale. Then for any $L>0$, it follows that
		\begin{flalign*}
			\PP \Big( \sup_{0 \le t \le T} Z_t \ge L \Big) \le\frac{ \E[Z_0]}{L}.
		\end{flalign*}
	\end{proposition}

    \begin{proposition}[Discrete Lenglart domination inequality \cite{1977,Mehri2021}]\label{DLDI}
    Let $\{Z_n\}_{n=0}^N$ be a nonnegative adapted process and
    $\{A_n\}_{n=0}^N$ an increasing predictable process with $A_0=0$.
     Suppose that $\mathbb E[Z_\tau]\le \mathbb E[A_\tau]$ for every stopping time $\tau\le N$. Then, for every $\theta\in(0,1)$,
     \begin{flalign*}
         \mathbb E\Big[ \Big(\max_{0\le n\le N}Z_n\Big)^\theta \Big]
    \le C_\theta\mathbb E[A_N^\theta].
     \end{flalign*}
    This is the discrete-time form of Lenglart's domination inequality;
    see \cite[Theorem~I and Corollay~II]{1977} and \cite[Lemma~2.2]{Mehri2021}.
    \end{proposition}
    
    Assume that $G\in C^1((A,B);\mathbb R)$ and $G(x)>0$ for $x\in(A,B)$. For a fixed constant $\sigma>0$, define the Lamperti-type transformation
	\begin{flalign*}
		\Phi(x)=\sigma\int_{d_1}^{x}\frac{1}{G(u)}\mathrm{d}u+d_{2},\quad x\in(A,B),
	\end{flalign*}
	where $d_{1}\in(A,B)$ and $d_{2}\in\RR$ are arbitrary constants. Then
    \begin{flalign*}
		\Phi'(x)=\frac{\sigma}{G(x)}\quad\text{and}\quad\Phi''(x)=-\frac{\sigma G'(x)}{G(x)^2}.
	\end{flalign*}
    Hence $\Phi$ is strictly increasing on $(A,B)$. Set
    \begin{equation*}
		a := \lim_{x \downarrow A} \Phi(x), \qquad b := \lim_{x \uparrow B} \Phi(x),
	\end{equation*}
    Then $\Phi:(A,B)\longrightarrow(a,b)$ is a bijection. Define
    \begin{equation*}
		\mu(w) = \sigma \left[ \frac{F(x)}{G(x)} - \frac{1}{2} G'(x) \right]_{x = \Phi^{-1}(w)},\quad w\in(a,b).
	\end{equation*}
Then the  Lamperti-type
transformation $W_t=\Phi(U_t)$ transforms the equation \eqref{sde} to 
    \begin{equation}\label{tsde}
		\mathrm{d}W_t = \mu(W_t)\,\mathrm{d}t + \sigma\,\mathrm{d}B_t, \quad W_0 = \Phi(U_0)\in(a,b).
	\end{equation}

    \begin{assumption}\label{as1}
		Assume that $\mu\in C^2_{loc}(a,b)$ and there exists $H>0$ such that 
		\begin{flalign*}
			\langle w_1-w_2, \mu(w_1)-\mu(w_2)\rangle\leq H|w_1-w_2|^2, \quad w_1,w_2\in(a,b).
		\end{flalign*}
	\end{assumption}
	\begin{assumption}\label{as2}
		If $a > -\infty$, assume that there exist $\bar{q}_a > 0$, $\varepsilon_a > 0$, and $C_a > 0$ such that
		\begin{equation*}
			(x - a)\mu(x) \ge \frac{\bar{q}_a + 1}{2}\sigma^2 + C_a, \qquad a < x \le a + \varepsilon_a.
		\end{equation*}
		If $b < \infty$, assume that there exist $\bar{p}_b > 0$, $\varepsilon_b > 0$, and $C_b > 0$ such that
		\begin{equation*}
			(b - x)\mu(x) \le -\frac{\bar{p}_b + 1}{2}\sigma^2 - C_b, \qquad b - \varepsilon_b \le x < b.
		\end{equation*}
	\end{assumption}
Fix $w^* \in  (a,b)$ and let $y=w^*$. Then from Assumption \ref{as1}, one obtains
	\begin{flalign}\label{xw*}
		(x-w^*)\mu(x)=&(x-w^*)(\mu(x)-\mu(w^*))+(x-w^*)\mu(w^*)\nonumber\\
		\leq&H|x-w^*|^2+\frac12|x-w^*|^2+\frac12|\mu(w^*)|^2\leq C(1+|x-w^*|^2).
	\end{flalign}
    Under the preceding assumptions on the SDE \eqref{tsde}, we obtain the following crucial properties for the solution $\{W_t\}_{t\in[0,T]}$. In what follows, we first establish the well-posedness of the exact solution.
	\begin{lemma}\label{Lm2.1}
		Let Assumptions \ref{as1} and \ref{as2} hold. Then SDE \eqref{tsde} admits a unique strong solution on $[0,T]$, and
		\begin{flalign*}
			\PP(W_t\in(a,b),\,0\leq t\leq T)=1.
		\end{flalign*}
	\end{lemma}
	\begin{proof}
		Since $\mu\in C^2_{loc}$, the SDE \eqref{tsde} admits a unique maximal local strong solution $\{W(t)\}_{[0, \tau_e)}$, where $\tau_e$ is the explosion time. For sufficiently large $n$, define
		\begin{flalign*}
			a_n = \begin{cases} a + n^{-1}, & a > -\infty, \\ -n, & a = -\infty, \end{cases}
			\qquad
			b_n = \begin{cases} b - n^{-1}, & b < \infty, \\ n, & b = \infty, \end{cases}
		\end{flalign*}
		and
		\begin{flalign*}
			\tau_n := \inf \{ t \in [0, \tau_e) : W_t \notin (a_n, b_n) \} \wedge \tau_e.
		\end{flalign*}
		Since $[a_n,b_n]\Subset(a,b)$, $\tau_n$ is increasing, and $\tau_n \uparrow \tau_e$ a.s. Fix $w^* \in (a, b)$ and define $V_0(x):=1+|x-w^*|^2$. By \eqref{xw*}, 
		\begin{flalign}\label{LV0}
			\mathcal LV_0(x) = 2(x-w^*)\mu(x)+\sigma^2 \le C(1+|x-w^*|^2)=CV_0(x).
		\end{flalign}
		If $a>-\infty$, choose $0<q\le\bar q_a$ and sufficiently small $\eta_a>0$ such that $2\eta_a<\varepsilon_a$ and $a+2\eta_a<b$. Let $\chi_a\in C^2(a,b)$ and $0\le\chi_a\le1$, which satisfies $\chi_a(x)=1$, $ a<x\le a+\eta_a$ and $\chi_a(x)=0$, $x\ge a+2\eta_a$. Define $V_a(x):=\chi_a(x)(x-a)^{-q}$. For $a<x\le a+\eta_a$, Assumption \ref{as2} gives
		\begin{flalign}\label{Va}
			\mathcal{L}V_a(x)=
			q(x-a)^{-q-2} \Big( \frac{q+1}{2}\sigma^2 - (x-a)\mu(x) \Big)\leq -C_aq(x-a)^{-q-2}\leq 0.
		\end{flalign} 
	    Since $[a+\eta_a,a+2\eta_a]\Subset(a,b)$, we have $\sup_{x\in[a+\eta_a,a+2\eta_a]} \mathcal LV_a(x)<\infty$. Moreover, $V_a(x)= 0$ for $x\ge a+2\eta_a$. Consequently, for $x\in(a,b)$, $\mathcal LV_a(x) \le C (1+V_a(x))$. Similarly, if $b<\infty$, choose $0<p\le\bar p_b$ and sufficiently small $\eta_b>0$, together with $\chi_b\in C^2(a,b)$ and $0\le\chi_b\le1$, such that $\chi_b(x)=1$, $b-\eta_b\le x<b$ and $\chi_b(x)=0$, $x\le b-2\eta_b$. Let $V_b(x):=\chi_b(x)(b-x)^{-p}$. Then, for $b-\eta_b\le x<b$,
		\begin{flalign}\label{Vb}
			\mathcal LV_b(x) \le p(b-x)^{-p-2} \Big( -\frac{\bar p_b+1}{2}\sigma^2-C_b +\frac{p+1}{2}\sigma^2 \Big)\le -C_bp(b-x)^{-p-2} \le 0.
		\end{flalign}
		As above, $\mathcal LV_b$ is bounded on $[b-2\eta_b,b-\eta_b]$ and $\mathcal LV_b(x) \le C\bigl(1+V_b(x)\bigr)$, $x\in(a,b)$.
		Now define $V(x) := V_0(x) + \mathbb I_{\{a>-\infty\}}V_a(x) + \mathbb I_{\{b<\infty\}}V_b(x)$. Since $V_0(x)\ge 1$, combining \eqref{LV0}, \eqref{Va}, and \eqref{Vb} gives 
	    \begin{flalign*}
	    	\mathcal LV(x)\le CV(x), \quad x\in(a,b).
	    \end{flalign*}
	    Applying It\^o's formula to $V(W_{t \wedge \tau_n})$ and taking the expectation yields:
		\begin{flalign*}
			\mathbb{E}[V(W_{t \wedge \tau_n})] =  \mathbb{E}[V(W_0)] + \mathbb{E} \left[ \int_0^{t \wedge \tau_n} \mathcal{L}V(W_s) \mathrm{d}s \right]
			\leq  \mathbb{E}[V(W_0)] + C\mathbb{E} \left[ \int_0^{t }V(W_{s\wedge \tau_n}) \mathrm{d}s\right].
		\end{flalign*}
		Gronwall's inequality yields
		\begin{flalign}\label{EVstop}
			\sup_{n\ge 1}\sup_{0\le t\le T} \mathbb E V(W_{t\wedge\tau_n}) \le C_T.
		\end{flalign}
		Define $c_n:=\min\{V(a_n),V(b_n)\}$. We claim that $c_n\rightarrow\infty$. Indeed, if $a>-\infty$, then for sufficiently large $n$ such that $n^{-1}<\eta_a$, which gives $\chi_a(a_n)=1$ and
		\begin{flalign*}
			V(a_n) \ge V_a(a_n) = (a_n-a)^{-q} = n^q \rightarrow\infty.
		\end{flalign*}
		 If $a=-\infty$, then
		 \begin{flalign*}
		 	V(a_n) \ge V_0(-n) = 1+\vert{}-n-w^*\vert{}^2 \rightarrow\infty.
		 \end{flalign*}
		 Similarly, if $b<\infty$, then for sufficiently large $n$,
		 \begin{flalign*}
		 	V(b_n) \ge V_b(b_n) = (b-b_n)^{-p} = n^p \rightarrow\infty,
		 \end{flalign*}
		 whereas if $b=\infty$, then 
		 \begin{flalign*}
		 	V(b_n) \ge V_0(n) = 1+\vert{}n-w^*\vert{}^2 \rightarrow\infty.
		 \end{flalign*}
		 Hence $c_n\to\infty$. Taking $t=T$ in \eqref{EVstop}, we obtain
		 \begin{flalign*}
		 	 c_n\mathbb P(\tau_n\le T) \le\mathbb E \left[ V(W_{\tau_n}) \mathbb I_{\{\tau_n\le T\}} \right]\leq  \mathbb E V(W_{T\wedge\tau_n})\leq C_T.
		 \end{flalign*}
		 Since $\tau_n\uparrow\tau_e$ a.s., we have $\{\tau_e\le T\}\subseteq\{\tau_n\le T\} $ for every $n$, and hence
		 \begin{flalign*}
		 	\mathbb P(\tau_e\le T) \leq \mathbb P(\tau_n\le T) \le \frac{C_T}{c_n} \rightarrow 0.
		 \end{flalign*}
		 Since $T>0$ is arbitrary, $\tau_e=\infty$ a.s. Thus the maximal local solution is global and remains in $(a,b)$ almost surely. Therefore,
		 \begin{flalign*}
		    \mathbb P\bigl(W_t\in(a,b),\ 0\le t\le T\bigr)=1.
		 \end{flalign*}
		 This completes the proof.
	\end{proof}
	
	With the global well-posedness of the transformed solution established, we proceed to derive its moment estimates.
	\begin{lemma}\label{Lm2.2}
		Let Assumptions \ref{as1} and \ref{as2} hold. For any $p >0$ and finite time $T > 0$, one has
		\begin{equation*}
			\mathbb{E} \Big[ \sup_{t \in [0,T]} |W_t|^p \Big] \leq C. 
		\end{equation*}
	\end{lemma}
    \begin{proof}
    For $p\geq 2$, the proof can be done in the same  way as  that of Lemma 3.2 in \cite{HZS2002} and is therefore omitted. For $0<p<2$, the result follows from the case $p=2$ and Jensen's inequality.
    \end{proof}
	
    Moment estimates provide control of the solution at infinity. In addition, the negative moment  estimates established in the following lemma are needed to characterize how closely the exact solution may approach the boundary. If $a>-\infty$, define $q_a^\dagger:=\bar q_a+\frac{2C_a}{\sigma^2}$. If $b<\infty$, define $p_b^\dagger:=\bar p_b+\frac{2C_b}{\sigma^2}$.
	\begin{lemma}\label{Lm2.3}
		Let Assumptions \ref{as1} and \ref{as2} hold. If $a>-\infty$, then for every $ 0<q<q_a^\dagger $, there exists a constant $C > 0$ such that
		\begin{flalign*}
			\mathbb E\Big[\sup_{0\le t\le T}(W_t-a)^{-q}
			\Big] \le C.
		\end{flalign*}
		If \(b<\infty\), then for every $0<p<p_b^\dagger$, there exists a constant $C > 0$ such that
		\begin{flalign*}
			\mathbb E\Big[\sup_{0\le t\le T}(b-W_t)^{-p}\Big]
			\le C.
		\end{flalign*}
	\end{lemma}
	\begin{proof}
		Let $a>-\infty$ and $0<q<q_a^\dagger$. Choose $q_0\in (q,q_a^\dagger)$ and set $d_a:= \frac{\bar q_a+1}{2}\sigma^2+C_a -\frac{q_0+1}{2}\sigma^2>0$. Let $\eta_a$ and $\chi_a$ be chosen as  in Lemma \ref{Lm2.1}. Define $\tilde{V}_a(x) = 1+ \chi_a(x)(x-a)^{-q_0}$. For $a<x\le a+\eta_a$, Assumption \ref{as2} gives
		\begin{flalign*}
			\mathcal L\widetilde V_a(x)\le -q_0d_a(x-a)^{-q_0-2} \le 0.
		\end{flalign*}
		On interval $[a+\eta_a,a+2\eta_a]$, $\mathcal L\widetilde V_a$ is bounded, while $\widetilde V_a(x)= 1$ for $x\ge a+2\eta_a$. Hence there exists $C_v>0$ such that
		\begin{flalign*}
			\mathcal L\widetilde V_a(x) \le C_v\widetilde V_a(x), \quad x\in(a,b).
		\end{flalign*}
		Recall $\tau_n$ and define $Q_t^{(n)}=e^{-C_v(t\wedge \tau_n)}\widetilde{V}_a(W_{t\wedge\tau_n})$. By It\^o's formula, 
		\begin{flalign*}
			Q_t^{(n)}=\widetilde{V}_a(W_0)+\int_{0}^{t\wedge \tau_n}e^{-C_{v}s}\big(\mathcal L \widetilde{V}_a(W_s)-C_{v}\widetilde{V}_a(W_s)\big)\mathrm{d}s+\sigma\int_{0}^{t\wedge \tau_n} e^{-C_{v}s}\widetilde{V}_a'(W_s)\mathrm{d}B_s.
		\end{flalign*}
		Consequently, for $0 \le s \le t \le T$,
		\begin{flalign*}
			\mathbb E\left[ Q_t^{(n)} \mid \mathcal F_s \right] = Q_s^{(n)} + \mathbb E\left[  \int_{s \wedge \tau_n}^{t \wedge \tau_n} e^{-C_vu} \bigl[\mathcal L\widetilde{V}_a(W_u) - C_v \widetilde{V}_a(W_u)\bigr] \mathrm du \Big| \mathcal F_s \right] \le Q_s^{(n)}.
		\end{flalign*}
		Thus $\{Q_t^{(n)}\}_{0 \le t \le T}$ is a non-negative supermartingale. By Proposition \ref{ville},
        \begin{flalign*}
            \mathbb P\Big( \sup_{0\le t\le T}Q_t^{(n)}>\lambda \Big)
    \le 1\wedge\frac{\mathbb E[Q_0^{(n)}]}{\lambda}, \quad \lambda>0.
        \end{flalign*}
        Set $K_0:=\mathbb E[Q_0^{(n)}] =\mathbb E[\widetilde V_a(W_0)]<\infty$, which is independent of $n$. Then, for every $\theta\in(0,1)$,
		\begin{flalign*}
			\mathbb E\Big[ \Big(\sup_{0\le t\le T}Q_t^{(n)}\Big)^\theta
        \Big] &= \theta\int_0^\infty \lambda^{\theta-1} \mathbb P\Big(
    \sup_{0\le t\le T}Q_t^{(n)}>\lambda \Big)\mathrm{d}\lambda\le \theta\int_0^{K_0}\lambda^{\theta-1}\,\mathrm{d}\lambda + \theta K_0 \int_{K_0}^{\infty}\lambda^{\theta-2}\,\mathrm{d}\lambda\\ &= K_0^\theta + \frac{\theta}{1-\theta}K_0^\theta= \frac{K_0^\theta}{1-\theta} \le C_\theta .
		\end{flalign*}
		where $C_\theta$ is independent of $n$. Since $\tau_n\uparrow\infty$ a.s. and
		$e^{C_v(t\wedge\tau_n)}\le e^{C_vT}$, Fatou's lemma yields
		\begin{flalign*}
			\mathbb E\Big[ \Big( \sup_{0\le t\le T} \widetilde{V}_a(W_t) \Big)^\theta
			\Big]  \le C.
		\end{flalign*}
		Moreover, $(x-a)^{-q_0} \le C\widetilde V_a(x)$, $ x\in(a,b)$. Taking $\theta=\frac{q}{q_0}$, we obtain
		\begin{flalign*}
			\mathbb E \Big[\sup_{0\le t\le T}(W_t-a)^{-q}\Big] = \mathbb E \Big[ \Big( \sup_{0\le t\le T} (W_t-a)^{-q_0} \Big)^{q/q_0} \Big]\le C\mathbb E \Big[ \Big( \sup_{0\le t\le T} \widetilde V_a(W_t) \Big)^{\theta} \Big] \le C.
		\end{flalign*}
	    The other estimate follows analogously by choosing $p<p_0<p_b^\dagger$ and using $\widetilde{V}_b=1+\chi_b(x)(b-x)^{-p_0}$.
	\end{proof}
    
	\section{The DPPEM scheme and strong convergence}
    \label{sec3}
	Having established the analytical estimates required to control the growth of the exact solution and its behavior near the boundary, we now construct the DPPEM approximation for the transformed SDE and establish its strong convergence.
	\subsection{The DPPEM scheme}
 
    For a given step size $h$, define the step-size-dependent interval $\mathbb{D}_h=[a_h, b_h]$, where
		\begin{flalign*}
			a_h = 
			\begin{cases} 
				a + h^{l_1}, & \text{if } a>-\infty,\\
				-h^{-m_a}, & \text{if } a=-\infty,
			\end{cases}
		\quad \text{and} \quad
			b_h = 
			\begin{cases} 
				b - h^{l_2}, & \text{if } b<\infty,\\
				h^{-m_b}, & \text{if } b=\infty,
			\end{cases}
		\end{flalign*}
	with $l_1,l_2,m_a,m_b>0$. Since $W_0\in(a,b)$, as well as $a_h\rightarrow a$ and $b_h\rightarrow b$ at infinite endpoints, there exists $0<h_0\leq 1$ such that $W_0\in\mathbb D_h\subset (a,b)$ for every $0<h\le h_0$. In particular, $\mathbb D_h$ is a nonempty closed convex interval for such $h$. Throughout this section, we assume $0<h\le h_0$ and set $W^0=W_0$.
    
    For $n=0,\ldots,N-1$, we first compute the explicit Euler predictor
\begin{equation}\label{preq}
    \widetilde W^{\,n+1}= W^n+h\mu(W^n)+\sigma\Delta B_n,
    \quad\Delta B_n=B(t_{n+1})-B(t_n).
\end{equation}
To enforce the domain constraint, we set
\begin{flalign*}
    W^{n+1}= \begin{cases}
     \widetilde W^{\,n+1}, & \quad \hbox{if $ \widetilde W^{\,n+1}\in \mathbb D_h$},\\
     a_h, & \quad \hbox{if $ \widetilde W^{\,n+1}< a_h$},\\
    b_h, & \quad \hbox{if $ \widetilde W^{\,n+1}> b_h$}.\\
\end{cases}
\end{flalign*}
This is the Euclidean projection of the predictor onto $\mathbb D_h$, namely $W^{n+1}=\Pi_{\mathbb D_h}(\widetilde W^{\,n+1})$, and can be evaluated explicitly. Since $W^0=W_0\in\mathbb D_h$, the scheme satisfies $W^n\in\mathbb D_h\subset(a,b)$ for $n=0,\ldots,N$,  pathwise. Therefore, the DPPEM approximation preserves the prescribed interval at every time step.

The EM step followed by a pointwise correction is related to the predictor--corrector approach of Cheng and Shen \cite{CPC2022,CPC2022.2} for positivity- and bound-preserving discretizations of parabolic equations. Their Lagrange multiplier formulations include variants in which the correction reduces to a cutoff operation. In the present scheme, the correction is an explicit projection onto the step-size-dependent interior interval $\mathbb D_h$, which keeps numerical values away from singular finite endpoints and imposes finite bounds on unbounded domains.
    
	\subsection{Strong convergence order of the DPPEM scheme}

    While the projection guarantees that the numerical approximation remains in the step-size-dependent interval $\mathbb D_h$, the strong convergence analysis requires quantitative control of the drift and its derivative near finite boundaries and at infinity whenever relevant. To distinguish the growth requirements arising from local consistency and stability, we use the subscripts \(I\) and \(S\), respectively. Specifically, \((\gamma_{I,a},\gamma_{I,b},\kappa_I)\) control the growth of \(\mu\) and \(L\mu\) in the local consistency estimates, whereas \((\gamma_{S,a},\gamma_{S,b},\kappa_S)\) control the growth of \(\mu'\) in the stability analysis. We therefore impose the following assumption.
	\begin{assumption}\label{as3}
		There exist nonnegative constants $\gamma_{I,a},\gamma_{I,b},\kappa_{I}$, $\gamma_{S,a},\gamma_{S,b},\kappa_{S}$, and a constant $C>0$ such that for all $x\in(a,b)$,
		\begin{flalign*}
			|\mu(x)|\vee |\mathcal L\mu(x)|\leq C(1+\mathbb{I}_{\{a > -\infty\}}(x-a)^{-\gamma_{I,a}}+\mathbb{I}_{\{b < \infty\}}(b-x)^{-\gamma_{I,b}}+\mathbb{I}_{\{a = -\infty \text{ or } b = \infty\}} |x|^{\kappa_I})
		\end{flalign*}
        and
        \begin{flalign*}
             |\mu'(x)|\leq C(1+\mathbb{I}_{\{a > -\infty\}}(x-a)^{-\gamma_{S,a}}+\mathbb{I}_{\{b < \infty\}}(b-x)^{-\gamma_{S,b}}+\mathbb{I}_{\{a = -\infty \text{ or } b = \infty\}} |x|^{\kappa_S}).
        \end{flalign*}
	\end{assumption}

    Combining the growth condition above with the moment and negative moment  estimates established in Lemmas~\ref{Lm2.2} and \ref{Lm2.3}, we obtain the integrability of the drift and its derivatives.
    \begin{lemma}\label{Integrability of the drift}
    	Let Assumptions \ref{as1}, \ref{as2}, and \ref{as3} hold. Let $Q>2$. Assume that
    	\begin{flalign}\label{qaqb}
    		\max \{Q\gamma_{I,a},Q\gamma_{S,a}\}<q_a^\dagger, \quad a>-\infty\quad\text{and}\quad
    		\max \{Q\gamma_{I,b},Q\gamma_{S,b}\}<p_b^\dagger, \quad b<\infty.
    	\end{flalign}
    	Then we have
    	\begin{flalign}\label{together}
    		\mathbb E\Big[ \sup_{0\le t\le T} \left( |\mu(W_t)|^{Q} + |\mu'(W_t)|^{Q} + |\mathcal L\mu(W_t)|^{Q} \right) \Big]\le C.
    	\end{flalign}
    \end{lemma}
    \begin{proof}
    	By Assumption \ref{as3}, we obtain
    	\begin{flalign*}
    		 |\mu(x)|^{Q} \vee |\mathcal L\mu(x)|^{Q} \le C\big( 1 + \mathbb I_{\{a>-\infty\}}(x-a)^{-Q\gamma_{I,a}} + \mathbb I_{\{b<\infty\}}(b-x)^{-Q\gamma_{I,b}}+ \mathbb I_{\{a=-\infty\text{ or }b=\infty\}} |x|^{Q\kappa_I} \big)
    	\end{flalign*}
        and
        \begin{flalign*}
    		 |\mu'(x)|^{Q} \le C\big( 1 + \mathbb I_{\{a>-\infty\}}(x-a)^{-Q\gamma_{S,a}} + \mathbb I_{\{b<\infty\}}(b-x)^{-Q\gamma_{S,b}}+ \mathbb I_{\{a=-\infty\text{ or }b=\infty\}} |x|^{Q\kappa_S} \big)
    	\end{flalign*}
    	If $a>-\infty$, then the inequality \eqref{qaqb} together with Lemma \ref{Lm2.3} gives 
    	\begin{flalign*}
    		\mathbb E \Big[\sup_{0\le t\le T} (W_t-a)^{-\max \{Q\gamma_{I,a},Q\gamma_{S,a}\}} \Big]\le C,
    	\end{flalign*}
    	and  similarly  if $b<\infty$,  we  get
    	\begin{flalign*}
    		\mathbb E \Big[\sup_{0\le t\le T} (b-W_t)^{-\max \{Q\gamma_{I,b},Q\gamma_{S,b}\}} \Big]\le C,
    	\end{flalign*}
    	If $a=-\infty$ or $b=\infty$, then Lemma \ref{Lm2.2} implies 
    	\begin{flalign*}
    		\mathbb E \Big[\sup_{0\le t\le T} |W_t|^{\max \{Q\kappa_I,Q\kappa_S\}}\Big] \le C.
    	\end{flalign*}
    	Combining the above estimates yields the assertion \eqref{together}.
    \end{proof}
    
    We next estimate the local consistency error of the Euler predictor.
    \begin{lemma}\label{RER}
    	Suppose that the assumptions of Lemma \ref{Integrability of the drift} hold. Then the exact solution satisfies
    	\begin{flalign*}
    		W(t_{n+1}) = W(t_{n})  + h\mu(W(t_{n}) ) + \sigma\Delta B_n + R_{n+1},
    	\end{flalign*}
    	where $R_{n+1} = \int_{t_n}^{t_{n+1}} \bigl(\mu(W_s) - \mu(W(t_{n}) )\bigr) \, \mathrm ds$. Moreover, 
    	\begin{flalign*}
    		||R_{n+1}||_{L^Q} \le Ch^{3/2} \quad \text{and}\quad ||\mathbb E_n R_{n+1}||_{L^Q} \le Ch^2.
    	\end{flalign*}
    \end{lemma}
    \begin{proof}
    	Applying It\^o's formula to $\mu(W_s)$ on $[t_n, s]$, we obtain
    	\begin{flalign*}
    		\mu(W_s) - \mu(W(t_n)) = \int_{t_n}^{s} \mathcal L\mu(W_u) \, \mathrm du + \sigma \int_{t_n}^{s} \mu'(W_u) \, \mathrm dB_u.
    	\end{flalign*}
    	Thus $R_{n+1} = R_{n+1}^{(1)} + R_{n+1}^{(2)} $, where
    	\begin{flalign*}
    		R_{n+1}^{(1)} = \int_{t_n}^{t_{n+1}} \int_{t_n}^{s} \mathcal L\mu(W_u) \, \mathrm du \, \mathrm ds \quad\text{and}\quad R_{n+1}^{(2)} = \sigma \int_{t_n}^{t_{n+1}} \int_{t_n}^{s} \mu'(W_u) \, \mathrm dB_u \, \mathrm ds.
    	\end{flalign*}
    	For $R_{n+1}^{(1)} $, Minkowski's inequality and Lemma \ref{Integrability of the drift} yield
    	\begin{flalign*}
    		||R_{n+1}^{(1)}||_{L^Q} \le \int_{t_n}^{t_{n+1}} \int_{t_n}^{s} ||\mathcal L\mu(W_u)||_{L^Q} \, \mathrm du \, \mathrm ds \le Ch^2.
    	\end{flalign*}
    	By stochastic Fubini's theorem,
    	\begin{flalign*}
    		R_{n+1}^{(2)} = \sigma \int_{t_n}^{t_{n+1}} (t_{n+1} - u)\mu'(W_u) \, \mathrm dB_u.
    	\end{flalign*}
    	Since $Q > 2$, we have $Q/2 > 1$. The BDG inequality gives
    	\begin{flalign*}
    		\Vert{}R_{n+1}^{(2)}\Vert{}_{L^Q}^2 &\le C \left\Vert{} \int_{t_n}^{t_{n+1}} (t_{n+1} - u)^2 \vert{}\mu'(W_u)\vert{}^2 \, \mathrm du \right\Vert{}_{L^{Q/2}} \le C \int_{t_n}^{t_{n+1}} (t_{n+1} - u)^2 \Vert{}\mu'(W_u)\Vert{}^2_{L^{Q}} \, \mathrm du \\ &= C \int_{t_n}^{t_{n+1}} (t_{n+1} - u)^2 \Vert{}\mu'(W_u)\Vert{}_{L^Q}^2 \, \mathrm du \le Ch^3.
    	\end{flalign*}
    	Therefore, $\Vert{}R_{n+1}^{(2)}\Vert{}_{L^Q} \le Ch^{3/2}$. Moreover, $R_{n+1}^{(2)}$ is an It\^o integral over $[t_n, t_{n+1}]$, and hence $\mathbb E_n R_{n+1}^{(2)} = 0$. Moreover, we have
    	\begin{flalign*}
    		||\mathbb E_n R_{n+1}||_{L^Q} = ||\mathbb E_n R_{n+1}^{(1)}||_{L^Q} \le ||R_{n+1}^{(1)}||_{L^Q} \le Ch^2.
    	\end{flalign*}
    	This completes the proof.
    \end{proof}
    To simplify the subsequent global error recursion, define
    \begin{flalign*}
    	l_n := 1 + h^{-2} |\mathbb E_n R_{n+1}| + h^{-3/2} \big( \mathbb E_n |R_{n+1}|^Q \big)^{1/Q}.
    \end{flalign*}
    \begin{lemma}\label{lQn}
    	Suppose that the conditions of Lemma \ref{Integrability of the drift} hold. Then we have 
    	\begin{flalign*}
    		 \sup_{0 \le n \le N-1} \mathbb E [l_n^Q] \le C. 
    	\end{flalign*}
    	\begin{proof}
    		Using the inequality $(x+y+z)^Q \le C_Q(x^Q+y^Q+z^Q)$, we obtain
    		\begin{flalign*}
    			\mathbb E [l_n^Q] \le C + C h^{-2Q} \mathbb E [\vert{}\mathbb E_n R_{n+1}\vert{}^Q] + C h^{-3Q/2} \mathbb E \Big[\Big(\big( \mathbb E_n \vert{}R_{n+1}\vert{}^Q \big)^{1/Q} \Big)^Q\Big].
    		\end{flalign*}
    		The last term simplifies exactly to
    		\begin{flalign*}
    			\mathbb E \Big[\Big(\big( \mathbb E_n \vert{}R_{n+1}\vert{}^Q \big)^{1/Q} \Big)^Q\Big] = \mathbb E[ \mathbb E_n \vert{}R_{n+1}\vert{}^Q] = \mathbb E \vert{}R_{n+1}\vert{}^Q.
    		\end{flalign*}
    		Therefore, by Lemma \ref{RER}, we derive
    		\begin{flalign*}
    			\mathbb E [l_n^Q] \le C + C h^{-2Q} h^{2Q} + C h^{-3Q/2} h^{3Q/2} \le C.
    		\end{flalign*}
    		This completes the proof.
    	\end{proof}
    \end{lemma}

    We now quantify the regularizing effect of the projection interval on the singular drift. Define
    \begin{flalign*}
        \delta:= \max \{ \mathbb I_{\{a > -\infty\}} l_1 \gamma_{S,a}, \, \mathbb I_{\{b < \infty\}} l_2 \gamma_{S,b}, \, \mathbb I_{\{a = -\infty\}} m_a \kappa_S, \, \mathbb I_{\{b = \infty\}} m_b \kappa_S\}.
    \end{flalign*}
    \begin{lemma}\label{hdelta}
    	 Suppose Assumption \ref{as3} holds. Then for all $x,y\in\mathbb D_h$, 
    	 \begin{flalign*}
    	 	|\mu(x) - \mu(y)|\le Ch^{-\delta} |x - y|.
    	 \end{flalign*}
    \end{lemma}
	\begin{proof}
		If $a > -\infty$, then $x - a \ge h^{l_1}$ for $x \in\mathbb D_h$, and hence $(x - a)^{-\gamma_{S,a}} \le h^{-l_1\gamma_{S,a}}$. Similarly, $(b - x)^{-\gamma_{S,b}} \le h^{-l_2\gamma_{S,b}} $ when $b < \infty$. If $a=-\infty$ or $b = \infty$, $|x|\le C ( 1 + \mathbb{I}_{\{a=-\infty\}}h^{-m_a} + \mathbb{I}_{\{b=\infty\}}h^{-m_b} )$. Hence, $|x|^{\kappa_S} \le Ch^{-\max\{\mathbb{I}_{\{a=-\infty\}}m_a \kappa_S, \, \mathbb{I}_{\{b=\infty\}}m_b \kappa_S\}}$. The mean-value theorem yields 
		\begin{flalign*}
		|\mu(x) - \mu(y)|\le Ch^{-\delta} |x - y|.
		\end{flalign*}
		This completes the proof.
	\end{proof}

    The preceding estimate provides an $h$-dependent local Lipschitz bound. The following lemma combines this bound with the one-sided Lipschitz condition.
    \begin{lemma}\label{true}
        Let Assumptions \ref{as1} and \ref{as3} hold and assume that
        $\delta<1$. Then for sufficiently small $h>0$, there exists a  constant $C>0$ independent of $h$, such that for any $x,y\in\mathbb D_h$,
        \begin{flalign*}
            |x-y+h(\mu(x)-\mu(y))|\leq (1+Ch)|x-y|.
        \end{flalign*}
    \end{lemma}
    \begin{proof}
        If $x=y$, the inequality above holds clearly. Let $a_h(x,y)=\frac{\mu(x)-\mu(y)}{x-y}$. One-sided Lipschitz condition gives $a_h(x,y)\leq H$. However, Lemma \ref{hdelta} means $|a_h(x,y)|\leq Ch^{-\delta}$, which gives $a_h(x,y)\geq -Ch^{-\delta}$. Since $\delta<1$, we can take sufficiently small $h$ such that $Ch^{1-\delta}\leq \frac12$. Hence $\frac12\leq 1+ha_h(x,y)\leq 1+Hh$, which means that
        \begin{flalign*}
            |x-y+h(\mu(x)-\mu(y))|=|1+ha_h(x,y)||x-y|\leq (1+Hh)|x-y|.
        \end{flalign*}
        This completes the proof.
    \end{proof}
    
    We first estimate the strong error before the exact solution leaves the projection interval. Define $e_n:=W(t_n)-W^n$ and $\sigma_h := \inf \{ n \in \{0, \ldots, N\}: W(t_n) \notin \mathbb D_h\}$, with $\inf \varnothing = N + 1$. 
	\begin{lemma}\label{EneQ}
		Suppose that the conditions of Lemma \ref{Integrability of the drift} hold. Let $Q > 2$ and  $\delta< 1 $. Then for sufficiently small $h$,
		\begin{flalign*}
			\mathbb E_n \left[ \vert{}e_{n+1}\vert{}^Q \mathbb I_{\{n+1 < \sigma_h\}} \right] \le \mathbb I_{\{n < \sigma_h\}} \left[ (1 + Ch)\vert{}e_n\vert{}^Q + Ch^{Q+1} (1 + l_n^Q) \right].
		\end{flalign*} 
	\end{lemma}
	\begin{proof}
		On the event $A_n:=\{n < \sigma_h\}$, $W(t_n)\in\mathbb D_h$, Lemma \ref{hdelta} gives 
		\begin{flalign*}
			|\mu(W(t_n))-\mu(W^n)| \le Ch^{-\delta} |e_n|.
		\end{flalign*}
		Define $\widehat W(t_{n+1}) := \Pi_{\mathbb D_h}(W(t_{n+1}))$ and $\widehat e_{n+1} := \widehat W(t_{n+1}) - W^{n+1}$. Since $\mathbb{D}_h$ is closed and convex, $\Pi_{\mathbb{D}_h}$ is non-expansive. Hence,
		\begin{flalign*}
			\vert{}\widehat e_{n+1}\vert{} = \vert{}\Pi_{\mathbb D_h}(W(t_{n+1})) - \Pi_{\mathbb D_h}(\widetilde W^{n+1})\vert{} \le \vert{}W(t_{n+1}) - \widetilde W^{n+1}\vert{}=|e_n + h( \mu(W(t_n))-\mu(W^n)) + R_{n+1}|.
		\end{flalign*}
		For any $Q > 2$, there exists $C_Q > 0$ such that
		\begin{flalign*}
			\vert{}x+y\vert{}^Q \le \vert{}x\vert{}^Q + Q\vert{}x\vert{}^{Q-2}xy + C_Q ( \vert{}x\vert{}^{Q-2}\vert{}y\vert{}^2 + \vert{}y\vert{}^Q).
		\end{flalign*}
		Applying this with $x = e_n+h( \mu(W(t_n))-\mu(W^n))$, $y = R_{n+1}$, and taking conditional expectation, we obtain
        \begin{flalign}\label{total}
        \mathbb E_n [ |\widehat e_{n+1}|^Q \mathbb I_{A_n} ] \le& \vert{}e_n + h( \mu(W(t_n))-\mu(W^n))\vert{}^Q \mathbb I_{A_n} \nonumber + Q\vert{}e_n + h( \mu(W(t_n))-\mu(W^n))\vert{}^{Q-1}  |\mathbb E_n [R_{n+1}]| \mathbb I_{A_n} \nonumber \\
        &+ C_Q\vert{}e_n + h( \mu(W(t_n))-\mu(W^n))\vert{}^{Q-2} \mathbb E_n [|R_{n+1}|^2] \mathbb I_{A_n} + C_Q \mathbb E_n [|R_{n+1}|^Q] \mathbb I_{A_n}.
        \end{flalign}
        By Lemma \ref{true}, we have
        \begin{flalign*}
            |e_n+h( \mu(W(t_n))-\mu(W^n))|\leq (1+Ch)|e_n|.
        \end{flalign*}
        Moreover, we derive
        \begin{flalign*}
            |e_n + h( \mu(W(t_n))-\mu(W^n))|^Q \mathbb I_{A_n}\leq (1+Ch)|e_n|^Q \mathbb I_{A_n}.
        \end{flalign*}
        Using the definition of $l_n$, we obtain
        \begin{flalign*}
            |e_n + h( \mu(W(t_n))-\mu(W^n))|^{Q-1}  |\mathbb E_n [R_{n+1}]| \mathbb I_{A_n}\leq Ch^2|e_n|^{Q-1}l_n\mathbb I_{A_n}\leq \varepsilon h|e_n|^Q\mathbb I_{A_n}+C_\varepsilon h^{Q+1}l^Q_n\mathbb I_{A_n}.
        \end{flalign*}
		By the definition of $l_n$, $(\mathbb E_n [|R_{n+1}|^Q])^{1/Q}\leq h^{3/2}l_n$. Therefore, it follows from Young's inequality that 
		\begin{flalign}\label{382}
			\vert{}e_n + h( \mu(W(t_n))-\mu(W^n))\vert{}^{Q-2} \mathbb E_n [|R_{n+1}|^2] \mathbb I_{A_n} \le C h^3|e_n|^{Q-2}l^2_n\leq \varepsilon h|e_n|^Q\mathbb I_{A_n} + C_\varepsilon h^{Q+1}l_n^Q\mathbb I_{A_n},
		\end{flalign}
		whereas 
        \begin{flalign*}
            \mathbb E_n [|R_{n+1}|^Q] \mathbb I_{A_n}\leq h^{3Q/2}l_n^Q\leq h^{Q+1}l^Q_n,
        \end{flalign*}
        where $\frac{3Q}{2}>Q+1$ is used. Combining these estimates into \eqref{total}, we derive
		\begin{flalign}\label{hat}
			\mathbb E_n [ |\widehat e_{n+1}|^Q\mathbb I_{A_n} ] \le \mathbb I_{A_n} [ (1+Ch)\vert{}e_n\vert{}^Q + Ch^{Q+1}(1+l_n^Q) ]. 
		\end{flalign}
		Since $A_{n+1}\subset A_n$ and $W(t_{n+1})\in\mathbb D_h$ on $A_{n+1}$, we have $\widehat W(t_{n+1})=W(t_{n+1})$ and $\widehat e_{n+1}=e_{n+1}$. Hence,
		\begin{flalign*}
			\vert{}e_{n+1}\vert{}^Q\mathbb I_{A_{n+1}} \le \vert{}\widehat e_{n+1}\vert{}^Q\mathbb I_{A_n}.
		\end{flalign*} 
		 Combining this with the inequality \eqref{hat} yields
		\begin{flalign*}
			\mathbb E_n [ \vert{}e_{n+1}\vert{}^Q\mathbb I_{A_{n+1}} ] \le \mathbb I_{A_n} [ (1+Ch)\vert{}e_n\vert{}^Q + Ch^{Q+1}(1+l_n^Q) ].
		\end{flalign*}
		This completes the proof.
	\end{proof}

    We next iterate the one-step estimate while preserving control of the
discrete-time maximum.
	\begin{lemma}\label{Emaxenr}
		Suppose that the assumptions of Lemma \ref{EneQ} hold. Then, for every $1 < r < Q$, 
		\begin{flalign}\label{Eenr}
			\mathbb E \Big[ \max_{\substack{0 \le n \le N \\ n < \sigma_h}} |e_n|^r \Big] \le C h^r.
		\end{flalign}
	\end{lemma}
	\begin{proof}
		Set $Y_n := |e_n|^Q \mathbb I_{A_n}$ and $\xi_n := C(1+l_n^Q)\mathbb I_{A_n}$. Lemma \ref{EneQ} gives
		\begin{flalign*}
			\mathbb E_n [Y_{n+1}] \le (1+Ch)Y_n + h^{Q+1}\xi_n.
		\end{flalign*}
		Define $\mathcal{G}_0 = 1$, $\mathcal{G}_{n+1} = \frac{\mathcal{G}_n}{1+Ch}$, and $\mathcal{S}_n := \mathcal{G}_n Y_n $. Multiplying this by $\mathcal{G}_{n+1}$, we obtain
		\begin{flalign*}
			\mathbb E_n [\mathcal{S}_{n+1}] \le \mathcal{S}_n + h^{Q+1}\mathcal{G}_{n+1}\xi_n.
		\end{flalign*}
		Define $J_0^* = 0$ and $ J_n^* := \sum_{j=0}^{n-1} h^{Q+1}\mathcal{G}_{j+1}\xi_j$. Then $\{\mathcal{S}_n - J_n^*\}_{n=0}^N$ is a supermartingale. Since $e_0 = 0$, we have $\mathcal{S}_0 = J_0^* = 0$. Then we have $\mathbb E[ \mathcal{S}_\nu] \le \mathbb E [J_\nu^*]$, where $\nu\leq N$ is arbitrary stopping index. Proposition \ref{DLDI} implies that for every $0 < \vartheta < 1$,
		\begin{flalign}\label{maxSn}
			\mathbb E \Big[ \Big( \max_{0 \le n \le N} \mathcal{S}_n \Big)^\vartheta \Big] \le C_\vartheta \mathbb E [(J_N^*)^\vartheta].
		\end{flalign}
	    Moreover, by $\mathcal{G}_n\leq 1$, we have
	    \begin{flalign*}
	    	J_N^*\leq Ch^{Q+1}\sum_{n=0}^{N-1}(1+l_n^Q)=Ch^Q\Big[h\sum_{n=0}^{N-1}(1+l_n^Q)\Big].
	    \end{flalign*}
	    Take $\vartheta=\frac rQ$. Consequently,
	    \begin{flalign*}
	    	\mathbb E [(J_N^*)^{r/Q}] \le C h^r \mathbb E \Bigg[\Big( h \sum_{n=0}^{N-1} (1+l_n^Q) \Big)^{r/Q}\Bigg].
	    \end{flalign*}
	    Since $\frac{r}{Q} < 1$, the function $x \mapsto x^{r/Q}$ is concave. Jensen's inequality gives
	    \begin{flalign}\label{EJN}
	    	\mathbb E [(J_N^*)^{r/Q} ]\le C h^r \Big( h \sum_{n=0}^{N-1} \mathbb E [1+l_n^Q] \Big)^{r/Q}.
	    \end{flalign}
	    By Lemma \ref{lQn} and $N h = T$, 
	    \begin{flalign*}
	    	h \sum_{n=0}^{N-1} \mathbb E (1+l_n^Q) \le C.
	    \end{flalign*}
	    Hence, $\mathbb E [(J_N^*)^{r/Q} ] \le C h^r$. On the other hand, $\mathcal{G}_n = (1+Ch)^{-n}$. Since $n \le N$ and $N h = T$, $\mathcal{G}_n \ge (1+Ch)^{-N} \ge e^{-CT} $. Therefore $Y_n \le e^{CT}\mathcal{S}_n$. Combining \eqref{maxSn}, \eqref{EJN}, and this yields
	    \begin{flalign*}
	    	\mathbb E \Big[\max_{0 \le n \le N} Y_n^{r/Q}\Big] \le C h^r.
	    \end{flalign*}
	    Since $Y_n^{r/Q} = \vert{}e_n\vert{}^r \mathbb I_{A_n}$, we obtain the assertion \eqref{Eenr}.
	\end{proof}

    It remains to estimate the probability that the exact solution leaves $\mathbb D_h$.
	\begin{lemma}\label{lem:exit-probability}
      Suppose that the conditions of Lemma \ref{Lm2.3} hold. Choose $P\geq 1$. If $a>-\infty$, choose $q_a^*\in(0,q_a^\dagger)$; if $b<\infty$, choose $p_b^*\in (0,p_b^\dagger)$. Define 
      \begin{flalign*}
          \Gamma_a=\left\{ \begin{array}{ll}l_1q_a^*,&a>-\infty,\\ m_aP,&a=-\infty, \end{array} \right.
          \quad
          \Gamma_b = \left\{ \begin{array}{ll} l_2p_b^*,&b<\infty,\\ m_bP,&b=\infty, \end{array} \right.
      \end{flalign*}
      and $\Gamma_X := \min \left\{\Gamma_a, \Gamma_b\right\}$. Then $\mathbb P(\sigma_h\le N) \le Ch^{\Gamma_X}$.
	  \end{lemma}
	 \begin{proof}
	 	If $a>-\infty$, then by Markov's inequality, we have
	 	\begin{flalign*}
	 		\mathbb P ( \min_{0\le n\le N}(W(t_n)-a)<h^{l_1} ) &\le \mathbb P ( \inf_{0\le t\le T}(W_t-a)<h^{l_1} )\le h^{l_1q_a^*} \mathbb E \Big[\sup_{0\le t\le T} (W_t-a)^{-q_a^*}\Big]\le Ch^{l_1q_a^*}.
	 	\end{flalign*}
	 	Similarly, if $b<\infty$,
	 	\begin{flalign*}
	 		\mathbb P ( \max_{0\le n\le N}W(t_n)>b-h^{l_2} ) \le Ch^{l_2p_b^*}.
	 	\end{flalign*}
	 	If $a=-\infty$, then
	 	\begin{flalign*}
	 		\mathbb P ( \min_{0\le n\le N}W(t_n)<-h^{-m_a} ) \le \mathbb P( \sup_{0\le t\le T}\vert{}W_t\vert{}>h^{-m_a} )\le h^{m_aP} \mathbb E \Big[\sup_{0\le t\le T}\vert{}W_t\vert{}^P\Big]\le Ch^{m_aP}.
	 	\end{flalign*}
	 	Likewise, if $b=\infty$,
	 	\begin{flalign*}
	 		\mathbb P ( \max_{0\le n\le N}W(t_n)>h^{-m_b} ) \le Ch^{m_bP}.
	 	\end{flalign*}
	 	Therefore, 
	 	\begin{flalign*}
	 		\mathbb P(\sigma_h\le N) \le C\Big( \mathbb I_{\{a>-\infty\}}h^{l_1q_a^*} + \mathbb I_{\{b<\infty\}}h^{l_2p_b^*}+ \mathbb I_{\{a=-\infty\}}h^{m_aP} + \mathbb I_{\{b=\infty\}}h^{m_bP} \Big).
	 	\end{flalign*}
	 	By the definition of $\Gamma_X$, for $0<h\le1$,
	 	\begin{flalign*}
	 	\mathbb P(\sigma_h\le N) \le Ch^{\Gamma_X}.
	 	\end{flalign*}
	 	This completes the proof.
	 \end{proof}
	
	We establish the global strong convergence of the DPPEM approximation.
    Let $r>1$ and suppose that the conditions of Lemma~\ref{Integrability of the drift} hold for some $Q>r$. We adopt the parameter choices and the definition of $\Gamma_X$ from \ref{lem:exit-probability}, with the additional requirement $P>r$ whenever at least one endpoint is infinite. Moreover, set $m^* := \max\{ \mathbb I_{\{a=-\infty\}} m_a,\, \mathbb I_{\{b=\infty\}} m_b \}$.
    
    \begin{theorem}\label{main}
    With the above notation and choices, assume that $\delta<1$. If both endpoints are finite, suppose that $\Gamma_X\ge r$. If at least one endpoint is infinite, suppose that
    \begin{flalign*}
        \Gamma_X\left(1-\frac{r}{P}\right)\ge r\quad\text{and}\quad \Gamma_X-m^*r\ge r.
    \end{flalign*}
    Then there exist $h_*\in(0,h_0]$ and $C>0$, independent of $h$, such that
\begin{flalign*}
    \mathbb E \left[ \max_{0\le n\le N} |W(t_n)-W^n|^r \right]
\le Ch^r, \quad 0<h\le h_*.
\end{flalign*}
Hence the DPPEM scheme converges strongly with order one.
    \end{theorem}
    \begin{proof}
    Choose $h_*\in(0,h_0]$ sufficiently small so that all small step size conditions required in Lemmas~\ref{true} and \ref{EneQ} hold. In the remainder of the proof, let $0<h\le h_*$.
    	Let $G_h:=\{\sigma_h=N+1\}$. We decompose
    	\begin{flalign}\label{EmaxGh}
    		\mathbb E \Big[\max_{0\le n\le N}\vert{}e_n\vert{}^r\Big] = \mathbb E \left[ \max_{0\le n\le N}\vert{}e_n\vert{}^r \mathbb I_{G_h} \right] + \mathbb E \left[ \max_{0\le n\le N}\vert{}e_n\vert{}^r \mathbb I_{G_h^c} \right].
    	\end{flalign}
    	On the event $G_h$, every $n<\sigma_h$, Lemma \ref{Emaxenr} gives 
    	\begin{flalign*}
    		\mathbb E \left[ \max_{0\le n\le N}\vert{}e_n\vert{}^r \mathbb I_{G_h} \right] \le Ch^r.
    	\end{flalign*}
    	If $a>-\infty$ and $b<\infty$, then $|e_n|\leq b-a$. We obtain
    	\begin{flalign*}
    		\mathbb E \left[ \max_{0\leq n\leq N}\vert{}e_n\vert{}^r\mathbb I_{G_h^c} \right] \le (b-a)^r\mathbb P(G_h^c) \le Ch^{\Gamma_X} \le Ch^r.
    	\end{flalign*}
    	If $a=-\infty$ or $b=\infty$, then for any $0\leq n\leq N$, we have $\vert{}W^n\vert{}\le C(1+h^{-m^*})$ and $\vert{}e_n\vert{}^r \le C(1+\vert{}W(t_n)\vert{})^r + Ch^{-m^*r}$. Applying H\"older's inequality with $P>r$ yields
    	\begin{flalign*}
    		\mathbb E \left[ \sup_{0\le t\le T} (1+\vert{}W_t\vert{})^r \mathbb I_{G_h^c} \right] \le \Big( \mathbb E \sup_{0\le t\le T} (1+\vert{}W_t\vert{})^P \Big)^{r/P} \mathbb P(G_h^c)^{1-r/P} \le Ch^{\Gamma_X(1-r/P)} \le Ch^r,
    	\end{flalign*}
    	and $h^{-m^*r}\mathbb P(G_h^c) \le Ch^{\Gamma_X-m^*r} \le Ch^r$. Therefore, combining these estimates with \eqref{EmaxGh} yields
    	\begin{flalign*}
    		\mathbb E \Big[\max_{0\le n\le N}\vert{}e_n\vert{}^r\Big]\leq Ch^r.
    	\end{flalign*}
    	This finishes the proof.
    \end{proof}

    As a consequence of the strong convergence result, we can further show that the numerical approximation inherits suitable moment and negative moment  bounds.
    \begin{lemma}\label{num-moment-inverse}
    	 Suppose the conditions of Theorem \ref{main} hold. Then
    	\begin{equation}\label{num-moment}
    		\mathbb E\Big[ \max_{0\le n\le N}|W^n|^r
    		\Big]\le C.
    	\end{equation}
    	Moreover, suppose that $a>-\infty$. Let $r_a>0$ satisfy $r_a<q_a^\dagger$. Let $q_a^*$ be as chosen in Theorem \ref{main}, and assume additionally that $ l_1r_a\Big(\frac1{q_a^*}+\frac1r\Big)<1$. Then
    	\begin{equation}\label{num-left-inverse}
    		\mathbb E\left[
    		\max_{0\le n\le N}(W^n-a)^{-r_a}
    		\right]\le C.
    	\end{equation}
    	Suppose that $b<\infty$. Let $r_b>0$ satisfy $r_b<p_b^\dagger$. Let $p_b^*$ be as chosen in Theorem \ref{main}, and assume additionally that $l_2r_b\Big(\frac1{p_b^*}+\frac1r\Big)<1$. Then
    	\begin{equation}\label{num-right-inverse}
    		\mathbb E\left[
    		\max_{0\le n\le N}(b-W^n)^{-r_b}
    		\right]\le C.
    	\end{equation}
    \end{lemma}
    \begin{proof}
    	We first prove \eqref{num-moment}. Since $W^n=W(t_n)-e_n$ and $r>1$,
    	\begin{flalign*}
    		|W^n|^r\le C_r\left(|W(t_n)|^r+|e_n|^r\right).
    	\end{flalign*}
    	Hence, by Lemma \ref{Lm2.2} and Theorem \ref{main},
    	\begin{align*}
    		\mathbb E\left[ \max_{0\le n\le N}|W^n|^r
    		\right]\le &C\mathbb E\left[
    		\max_{0\le n\le N}|W(t_n)|^r \right]
    		+ C\mathbb E\left[ \max_{0\le n\le N}|e_n|^r \right]\\
    		\le &C\mathbb E\left[ \sup_{0\le t\le T}|W_t|^r \right]
    		+ Ch^r\leq C.
    	\end{align*}
    	We next prove the negative moment  estimate. When $a>-\infty$, define $D_{a,h}:=\min_{0\le n\le N}\bigl(W(t_n)-a\bigr)$ and $ E_h:=\max_{0\le n\le N}|e_n|$.
    	Since the exact solution remains in $(a,b)$ almost surely, $	D_{a,h}>0$ a.s. Introduce the event $\Omega_{a,h}:= \{ E_h\le\frac12D_{a,h} \}$.
    	Then
    	\begin{align*}
    		\mathbb E\left[
    		\max_{0\le n\le N}(W^n-a)^{-r_a}
    		\right]
    		={}
    		\mathbb E\left[
    		\max_{0\le n\le N}(W^n-a)^{-r_a}
    		\mathbb I_{\Omega_{a,h}}
    		\right]
    		+
    		\mathbb E\left[
    		\max_{0\le n\le N}(W^n-a)^{-r_a}
    		\mathbb I_{\Omega_{a,h}^c}
    		\right].
    	\end{align*}
    	On $\Omega_{a,h}$, for every $n$,
    	\begin{align*}
    		W^n-a = W(t_n)-a-e_n \ge W(t_n)-a-|e_n| \ge
    		W(t_n)-a-E_h \ge \frac12\bigl(W(t_n)-a\bigr).
    	\end{align*}
    	Therefore, $(W^n-a)^{-r_a} \le 2^{r_a}(W(t_n)-a)^{-r_a}$
    	and hence, since $r_a<q_a^\dagger$, Lemma \ref{Lm2.3} gives
    	\begin{align}\label{left-good}
    		\mathbb E\left[
    		\max_{0\le n\le N}(W^n-a)^{-r_a}
    		\mathbb I_{\Omega_{a,h}}
    		\right] \le 2^{r_a} \mathbb E\Big[
    		\sup_{0\le t\le T}(W_t-a)^{-r_a}\Big] \le C.
    	\end{align}
    	On the other hand, since $W^n\in\mathbb D_h$ and $a_h=a+h^{l_1}$, we have $W^n-a\ge h^{l_1}$.
    	Thus, $\max_{0\le n\le N}(W^n-a)^{-r_a}
    	\le h^{-l_1r_a}$, and consequently
    	\begin{equation}\label{left-bad}
    		\mathbb E\left[ \max_{0\le n\le N}(W^n-a)^{-r_a}
    		\mathbb I_{\Omega_{a,h}^c}
    		\right] \le h^{-l_1r_a} \mathbb P(\Omega_{a,h}^c).
    	\end{equation}
    	By $ l_1r_a\left(\frac1{q_a^*}+\frac1r\right)<1$, $\frac{l_1r_a}{q_a^*} < 1-\frac{l_1r_a}{r}$.
    	Therefore, there exists $\eta_a\in(0,1)$ such that $\frac{l_1r_a}{q_a^*} < \eta_a < 1-\frac{l_1r_a}{r}$. Since $\Omega_{a,h}^c = \{ E_h>\frac12D_{a,h}
    	\}$, we have $\Omega_{a,h}^c
    	\subset \{D_{a,h}<h^{\eta_a}\}
    	\cup \{ E_h>\frac12h^{\eta_a}\}$.
    	Therefore, $\mathbb P(\Omega_{a,h}^c) \le
    	\mathbb P(D_{a,h}<h^{\eta_a}) + \mathbb P(
    	E_h>\frac12h^{\eta_a} )$. For the first term, Markov's inequality yields
    	\begin{align*}
    		\mathbb P(D_{a,h}<h^{\eta_a}) = \mathbb P(
    		D_{a,h}^{-q_a^*} 	> h^{-\eta_aq_a^*} ) \le
    		h^{\eta_aq_a^*} \mathbb E [D_{a,h}^{-q_a^*}] \le
    		h^{\eta_aq_a^*} \mathbb E[ \sup_{0\le t\le T}(W_t-a)^{-q_a^*} ]\le Ch^{\eta_aq_a^*},
    	\end{align*}
    	where Lemma \ref{Lm2.3} is used. For the second term, Markov's inequality and Theorem \ref{main} give
    	\begin{align*}
    		\mathbb P (
    		E_h>\frac12h^{\eta_a} ) \le
    		2^rh^{-\eta_ar}\mathbb E E_h^r &\le
    		Ch^{-\eta_ar}h^r
    		=
    		Ch^{r(1-\eta_a)}.
    	\end{align*}
    	Hence, $\mathbb P(\Omega_{a,h}^c) \le
    	Ch^{\eta_aq_a^*} + Ch^{r(1-\eta_a)}$.
    	Substituting this into \eqref{left-bad}, we obtain
    	\begin{align*}
    		\mathbb E\left[
    		\max_{0\le n\le N}(W^n-a)^{-r_a}
    		\mathbb I_{\Omega_{a,h}^c}
    		\right] \le
    		C h^{\eta_aq_a^*-l_1r_a} +
    		C h^{r(1-\eta_a)-l_1r_a}.
    	\end{align*}
    	By $\frac{l_1r_a}{q_a^*} < \eta_a < 1-\frac{l_1r_a}{r}$, $\eta_aq_a^*-l_1r_a>0$, and $r(1-\eta_a)-l_1r_a>0$.
    	Thus,
    	\begin{flalign*}
    		\mathbb E\left[ \max_{0\le n\le N}(W^n-a)^{-r_a}
    		\mathbb I_{\Omega_{a,h}^c} \right] \le C.
    	\end{flalign*}
    	Combining this estimate with \eqref{left-good} proves \eqref{num-left-inverse}. 
        
        Now we prove the right-endpoint estimate \eqref{num-right-inverse} by the same argument used for the left-endpoint estimate \eqref{num-left-inverse}. For $b<\infty$, define $D_{b,h}:=\min_{0\le n\le N}\bigl(b-W(t_n)\bigr)$ and $\Omega_{b,h} := \{ E_h\le\frac12D_{b,h} \}$.
    	On $\Omega_{b,h}$, $b-W^n \ge \frac12\bigl(b-W(t_n)\bigr)$,
    	and hence
    	\begin{flalign*}
    		\mathbb E\left[ \max_{0\le n\le N}(b-W^n)^{-r_b}
    		\mathbb I_{\Omega_{b,h}} \right] \le C
    	\end{flalign*}
    	by Lemma \ref{Lm2.3}. On $\Omega_{b,h}^c$, since $W^n\le b_h=b-h^{l_2}$, we have $b-W^n\ge h^{l_2}$, and therefore $\max_{0\le n\le N}(b-W^n)^{-r_b}
    	\le h^{-l_2r_b}$. The inequality $l_2r_b\Big(\frac1{p_b^*}+\frac1r\Big)<1$ implies that we can choose $\eta_b\in(0,1)$ such that $\frac{l_2r_b}{p_b^*} < \eta_b < 1-\frac{l_2r_b}{r}$. Similarly, we obtain
    	\begin{flalign*}
    		\mathbb P(\Omega_{b,h}^c) \le Ch^{\eta_bp_b^*}
    		+ Ch^{r(1-\eta_b)},
    	\end{flalign*}
    	and consequently
    	\begin{align*}
    		\mathbb E\left[
    		\max_{0\le n\le N}(b-W^n)^{-r_b}
    		\mathbb I_{\Omega_{b,h}^c}
    		\right]\le C h^{\eta_bp_b^*-l_2r_b} +
    		C h^{r(1-\eta_b)-l_2r_b} \le C.
    	\end{align*}
    	Therefore,
    	\begin{flalign*}
    		\mathbb E\left[ \max_{0\le n\le N}(b-W^n)^{-r_b}
    		\right]\le C.
    	\end{flalign*}
    	This completes the proof.
    \end{proof}
    
    \begin{remark}
    	Since $q_a^*$ can be chosen arbitrarily close to $q_a^\dagger$, condition $ l_1r_a\Big(\frac1{q_a^*}+\frac1r\Big)<1$ is satisfied for every $r_a>0$ such that $r_a<
    	\min \{ q_a^\dagger,\,
    	\frac{r q_a^\dagger}{l_1(r+q_a^\dagger)} \}$.
    	Similarly, the right-endpoint negative moment  estimate holds for every $r_b>0$ satisfying $r_b<
    	\min \{ p_b^\dagger,\,
    	\frac{r p_b^\dagger}{l_2(r+p_b^\dagger)} \}$.
    \end{remark}
	
    \section{Application to representative models}
    \label{sec4}
    In this section, we apply the convergence results of Sections \ref{sec3} to representative domain-constrained SDEs. For each model, we proceed in three steps:
    \begin{enumerate}[(\roman*)]
    \item Transform the equation and verify domain invariance (check Assumptions \ref{as1} and \ref{as2});
    \item Verify the strong convergence conditions (Assumption \ref{as3} and the conditions in Theorem \ref{main});
    \item Transform back and state the convergence result.
    \end{enumerate}

    \subsection{CIR process}
    Consider
\begin{flalign*}
 \mathrm{d}X_t=\kappa(\theta-X_t)\,\mathrm{d}t+\beta\sqrt{X_t}\,\mathrm{d}B_t, \quad X_0>0,\quad \kappa,\theta,\beta>0.
\end{flalign*}

(\romannumeral1)  Define
$\Phi:(0,\infty)\rightarrow(0,\infty)$ by $\Phi(x)=2\sqrt{x}$, with $\Phi^{-1}(w)=\frac{w^2}{4}$. Then the transformed equation is
\begin{flalign*}
 \mathrm{d}W_t
 =
 \mu(W_t)\,\mathrm{d}t+\beta\,\mathrm{d}B_t,
 \quad W_0=\Phi(X_0),
\end{flalign*}
where
\begin{flalign*}
 \mu(w)=2(\kappa\theta-\frac{\beta^2}{4})w^{-1}
 -\frac{\kappa}{2}w,
 \quad w>0.
\end{flalign*}
If $\frac{\kappa\theta}{\beta^2}>\frac12$, then $\mu'(w)<0$ for $w>0$. So Assumption \ref{as1} holds with $H=0$. Moreover,
\begin{flalign*}
 w\mu(w) =  2\kappa\theta-\frac{\beta^2}{2}  -\frac{\kappa}{2}w^2.
\end{flalign*}
Hence, for every
\begin{flalign*}
 0<q_a^\dagger<
 \frac{4\kappa\theta}{\beta^2}-2,
\end{flalign*}
Assumption \ref{as2} is  satisfied at the left endpoint $a=0$.
Therefore, Lemma \ref{Lm2.1} yields
\begin{flalign*}
 \mathbb P(W_t>0,\ 0\le t\le T)=1,
\end{flalign*}
and consequently $X_t>0$ almost surely.

(\romannumeral2) 
For the transformed drift,
\begin{flalign*}
 |\mu(w)|&\le C(1+w^{-1}+w),\\
 |\mu'(w)|&\le C(1+w^{-2}),\\
 |\mathcal L\mu(w)|&\le C(1+w^{-3}+w).
\end{flalign*}
Thus Assumption \ref{as3} holds with $\gamma_{I,a}=3$, $\gamma_{S,a}=2$, $\kappa_I=1$, $\kappa_S=0$. Accordingly, Lemma \ref{Integrability of the drift} requires $3Q<q_a^\dagger$ for some $Q>2$. Let $r\ge1$ be the desired moment order for the original CIR equation.
Since $\Phi^{-1}(w)=w^2/4$, we apply Theorem \ref{main} at the transformed
moment order $R=2r$. Assume $\frac{4\kappa\theta}{\beta^2}-2>6r $. Then we may choose
\begin{flalign*}
 R<Q< \frac{4\kappa\theta}{3\beta^2}-\frac23,
 \quad \text{and} \quad 3Q<q_a^\dagger< \frac{4\kappa\theta}{\beta^2}-2.
\end{flalign*}
Next choose $q_a^*\in(0,q_a^\dagger)$ and $l_1$ such that $\frac{R}{q_a^*}<l_1<\frac12$. Then $l_1q_a^*>R$ and $\delta=l_1\gamma_{S,a}=2l_1<1$.  Fix $R<c<l_1q_a^*$. Since Lemma \ref{Lm2.2} provides moments of arbitrary
positive order, choose $P>R$ sufficiently large such that $c(1-\frac{R}{P})\ge R$, and set $m_b=c/P$. Then $\Gamma_X=\min\{l_1q_a^*,m_bP\}=c$, and therefore
\begin{flalign*}
 \Gamma_X(1-\frac{R}{P})\ge R
 \qquad \text{and}\quad \Gamma_X-m_bR\ge R.
\end{flalign*}
Hence all the conditions of Theorem \ref{main} are satisfied at moment
order $R=2r$.

(\romannumeral3) Theorem \ref{main} gives
\begin{flalign*}
\E\!\left[\max_{0\le n\le N}|W(t_n)-W_n|^{2r}\right]\le Ch^{2r}. 
\end{flalign*}
Since
\begin{flalign*}
 |X(t_n)-X^n|^r=4^{-r}|W(t_n)-W_n|^r|W(t_n)+W_n|^r, 
\end{flalign*}
H\"older's inequality, together with Lemmas \ref{Lm2.2} and \ref{num-moment-inverse}, yields
\begin{flalign*}
\E\!\left[\max_{0\le n\le N}|X(t_n)-X^n|^r\right]\le Ch^r. 
\end{flalign*}
Condition $\frac{4\kappa\theta}{\beta^2}-2>6r$ is equivalent to
$ \frac{\kappa\theta}{\beta^2}>2$ and $1\le r<\frac16(\frac{4\kappa\theta}{\beta^2}-2)$.
Now we can apply Theorem \ref{main} to obtain  the following Proposition. 

\begin{proposition}\label{pro41}
Suppose $\frac{\kappa\theta}{\beta^2}>2$. For each $1\le r<\frac16(\frac{4\kappa\theta}{\beta^2}-2)$, choose the projection parameters as specified in step (ii) above. Then the DPPEM approximation $X^n=\Phi^{-1}(W_n)$ satisfies, for all sufficiently small $h>0$,
\begin{flalign*}
\E\!\left[\max_{0\le n\le N}|X(t_n)-X^n|^r\right]\le Ch^r, 
\end{flalign*}
where $C$ is independent of $h$. Hence the scheme converges strongly with order one for the original CIR model.
\end{proposition}

    \subsection{Heston $3/2$ process}

Consider
\begin{flalign*}
\mathrm{d}X_t=a_1X_t(a_2-X_t)\,\mathrm{d}t+a_3X_t^{3/2}\,\mathrm{d}B_t,
\quad X_0>0,\quad a_1,a_2,a_3>0.
\end{flalign*}

(\romannumeral1) Define $\Phi:(0,\infty)\to(-\infty,0)$ by $\Phi(x)=-2x^{-1/2}$, with $\Phi^{-1}(w)=4w^{-2}$. The associated transformed equation on $(-\infty,0)$ is
\begin{flalign*}
\mathrm{d}W_t=\mu(W_t)\,\mathrm{d}t+a_3\,\mathrm{d}B_t,\quad W_0=\Phi(X_0),
\end{flalign*}
where $\mu(w)=2(a_1+\frac{3a_3^2}{4})w^{-1}-\frac{a_1a_2}{2}w$. We have $\mu'(w)<0$, so Assumption \ref{as1} holds. For the finite right endpoint $b=0$ and every $0<p_b^\dagger < \frac{4a_1}{a_3^2}+2$, Assumption \ref{as2} can be satisfied at $b=0$.
Since $\frac{4a_1}{a_3^2}+2>0$, no additional
parameter restriction is required for domain invariance.
Therefore, Lemma \ref{Lm2.1} gives
\begin{flalign*}
 \mathbb P(W_t<0,\ 0\le t\le T)=1, 
\end{flalign*}
and therefore $X_t>0$ almost surely.

(\romannumeral2) Assumption \ref{as3} holds with $\gamma_{I,b}=3$, $\gamma_{S,b}=2$, $\kappa_I=1$, and $\kappa_S=0$. Accordingly, Lemma \ref{Integrability of the drift} requires  $3Q<p_b^\dagger$ for some $Q>2$. Let $r\ge1$ be the desired moment order for the original Heston $3/2$
equation. Since $\Phi^{-1}(w)=4w^{-2}$ is singular at $b=0$, both a
transformed strong error estimate and suitable negative moment  estimates are
required. Set $ R:=\frac{7r}{3}$ and $S:=\frac{21r}{4}$. Assume $\frac{4a_1}{a_3^2}+2 >7r$. Then one can choose $R<Q<\frac{4a_1}{3a_3^2}+\frac23$ and $3Q<p_b^\dagger<\frac{4a_1}{a_3^2}+2$. Since $Q>R$, we have $p_b^\dagger>3Q>7r$. Hence one may choose $7r<p_b^*<p_b^\dagger$. With the above choices of $R$ and $S$,
\begin{flalign*}
(p_b^*)^2-S(R+p_b^*) = (p_b^*-7r) \Big(p_b^*+\frac{7r}{4}\Big)>0.
\end{flalign*}
Therefore, $\frac{R}{p_b^*} < l_2 < \frac{1} {S(\frac1{p_b^*}+\frac1R)}$. Moreover, $\frac{1}{S(\frac1{p_b^*}+\frac1R)}<\frac{R}{S}=\frac49<\frac12$. Thus one can choose $l_2$ such that
\begin{flalign*}
\frac{R}{p_b^*}<l_2<\frac{1}{S(\frac1{p_b^*}+\frac1R)}<\frac12.
\end{flalign*}
Consequently, $l_2p_b^*>R$, $\delta=l_2\gamma_{S,b}=2l_2<1$, and $l_2S
(\frac1{p_b^*}+\frac1R)<1$, which is precisely the condition required in
Lemma \ref{num-moment-inverse} for the numerical inverse moment of order $S$. Fix $R<c<l_2p_b^* $. Since Lemma \ref{Lm2.2} provides moments of arbitrary positive order for $W_t$, choose $P>R$ sufficiently large such that $c(1-\frac{R}{P})\ge R$, and set $m_a=c/P$. Then $\Gamma_X=\min\{l_2p_b^*,m_aP\}=c$, and therefore
\begin{flalign*}
\Gamma_X\Big(1-\frac{R}{P}\Big)\ge R \quad\text{and}\quad
\Gamma_X-m_aR\ge R. 
\end{flalign*}
Thus all conditions of Theorem \ref{main} are satisfied at the transformed moment order $R=\frac{7r}{3}$.

(\romannumeral3) Theorem \ref{main} gives
\begin{flalign*} 
 \E\left[ \max_{0\le n\le N}|W(t_n)-W_n|^{R} \right] \le Ch^{R}.
\end{flalign*}
Since $S=\frac{21r}{4}<7r<p_b^\dagger$, Lemma \ref{Lm2.3} gives
\begin{flalign*}
\E\Big[\sup_{0\le t\le T}|W_t|^{-S}\Big]\le C,
\end{flalign*}
while the condition verified in (\romannumeral2) and
Lemma \ref{num-moment-inverse} yield
\begin{flalign*}
\E\left[
\max_{0\le n\le N}|W_n|^{-S}
\right]\le C.
\end{flalign*}
Since $|\frac4{w^2}-\frac4{v^2}| \le 8|w-v|(|w|^{-3}+|v|^{-3})$ for $w,v<0$, using H\"older's inequality gives
\begin{flalign*}
\E\left[\max_{0\le n\le N}|X(t_n)-X^n|^r\right]\le &C
\left(\E\left[\max_{0\le n\le N}|W(t_n)-W_n|^{R}\right]\right)^{\frac37}
\Big[\Big(\E\sup_{t\le T}|W_t|^{-S}\Big)^{\frac47}+
\Big(\E\max_{0\le n\le N}|W_n|^{-S}\Big)^{\frac47}\Big]\\
\le& Ch^r. 
\end{flalign*}
Thus $\frac{4a_1}{a_3^2}+2>7r$ or equivalently, $r<
\frac17 ( \frac{4a_1}{a_3^2}+2)$ is a sufficient condition.  
In particular, the admissible range contains $r=1$ whenever $\frac{a_1}{a_3^2}>\frac54$.
Then Theorem  \ref{main} applied to this model becomes  

\begin{proposition}\label{pro42}
Suppose $\frac{a_1}{a_3^2}>\frac54$. For each $1\le r<\frac17(\frac{4a_1}{a_3^2}+2)$, choose the projection parameters as specified in step (ii) above. Then the DPPEM approximation satisfies, for all sufficiently small $h>0$,
\begin{flalign*}
\E\left[\max_{0\le n\le N}|X(t_n)-X^n|^r\right]\le Ch^r. 
\end{flalign*}
Hence, the DPPEM scheme converges strongly with order one for the original Heston 3/2 model.
\end{proposition}

\subsection{CEV process}
Consider
\begin{flalign}\label{eq:cev}
\mathrm{d}X_t=\kappa(\theta-X_t)\,\mathrm{d}t+\beta X_t^\alpha\,\mathrm{d}B_t,
\quad X_0>0,\quad \kappa,\theta,\beta>0,\quad \frac12<\alpha<1.
\end{flalign}

(\romannumeral1) Put $\nu:=\frac1{1-\alpha}>2$ and define $\Phi:(0,\infty)\to(0,\infty)$ by $\Phi(x)=\frac{x^{1-\alpha}}{1-\alpha}$, and $\Phi^{-1}(w)=((1-\alpha)w)^\nu$. Then the transformed equation is
\begin{flalign*}
\mathrm{d}W_t
=\mu(W_t)\,\mathrm{d}t+\beta\,\mathrm{d}B_t,
\quad W_0=\Phi(X_0),
\end{flalign*}
where
\begin{flalign*}
\mu(w)
=&\,
\kappa\theta\bigl((1-\alpha)w\bigr)^{-\frac{\alpha}{1-\alpha}}
-\kappa(1-\alpha)w
-\frac{\alpha\beta^2}{2(1-\alpha)w}.
\end{flalign*}
As $w\downarrow0$, $w\mu(w)\to+\infty$ since $1-2\alpha<0$. Hence Assumption \ref{as2} can be satisfied at $a=0$ with arbitrarily large $q_a^\dagger$, and Lemma \ref{Lm2.3} consequently provides inverse moments of arbitrary positive order. Moreover, $\mu'(w)\to-\infty$ as $w\downarrow0$ and $\mu'(w)\to-\kappa(1-\alpha)$ as $w\to\infty$. Since $\mu'$ is continuous on $(0,\infty)$, it is bounded from above.
Thus Assumption \ref{as1} holds. Lemma \ref{Lm2.1} therefore yields a
unique global transformed solution in $(0,\infty)$, and
$X_t=\Phi^{-1}(W_t)>0$ almost surely.  

(\romannumeral2) The exponents are $\gamma_{S,a}=\frac1{1-\alpha}$, $\gamma_{I,a}=\frac{1+\alpha}{1-\alpha}$, $\kappa_I=1$, and $\kappa_S=0$. Fix $r\ge1$ and set $R:=2\nu r$. Choose $Q>\max\{R,2\}$ and $l_1\in(0,\frac1\nu)$. Since $q_a^\dagger$ can be chosen arbitrarily large, we may choose it so that $q_a^\dagger>\max\{Q\gamma_{I,a},\,\frac{R}{l_1}\}$. Then choose $\frac{R}{l_1}<q_a^*<q_a^\dagger$. It follows that $l_1q_a^*>R$ and $\delta =l_1\gamma_{S,a} =l_1\nu<1$. Fix $c\in(R,l_1q_a^*)$. Since Lemma \ref{Lm2.2} provides
moments of arbitrary positive order for $W_t$, choose $P>R$
sufficiently large such that $c\left(1-\frac{R}{P}\right)\ge R $, and set $m_b:=\frac{c}{P}$.
Then
\begin{flalign*}
\Gamma_X
=
\min\{l_1q_a^*,m_bP\}
=c,
\end{flalign*}
and hence
\begin{flalign*}
\Gamma_X\Big(1-\frac{R}{P}\Big)\ge R \quad\text{and}\quad
\Gamma_X-m_bR\ge R. 
\end{flalign*}
Thus all conditions of Theorem \ref{main} are satisfied at the
transformed moment order $R=2\nu r$.

(\romannumeral3) Theorem \ref{main} gives
\begin{flalign*}
\E\left[
\max_{0\le n\le N}|W(t_n)-W_n|^R
\right]
\le Ch^R.
\end{flalign*}
Since $R=2\nu r>2r$, it follows that
\begin{flalign*}
\E\left[
\max_{0\le n\le N}|W(t_n)-W_n|^{2r}
\right]
\le Ch^{2r}.
\end{flalign*}
Moreover, $2r(\nu-1)<R$, so Lemmas \ref{Lm2.2} and \ref{num-moment-inverse} yield
\begin{flalign*}
\E\left[
\max_{0\le n\le N}
\{1+|W(t_n)|+|W_n|\}^{2r(\nu-1)}
\right]\le C.
\end{flalign*}
Because $|\Phi^{-1}(w)-\Phi^{-1}(v)| \le C|w-v|\bigl(w^{\nu-1}+v^{\nu-1}\bigr)$ for $ w,v>0$, H\"older's inequality gives
\begin{flalign*}
 \E\Big[\max_{0\leq n\leq N}|X(t_n)-X^n|^r\Big]\le &C
\left(\E\Big[\max_{0\leq n\leq N}|W(t_n)-W_n|^{2r}\Big]\right)^{\frac12}
\left(\E\max_{0\leq n\leq N}\{1+|W(t_n)|+|W_n|\}^{2r(\nu-1)}\right)^{\frac12}\\
\leq &Ch^r.
\end{flalign*}

\begin{proposition}
Let $\kappa,\theta,\beta>0$ and $\frac12<\alpha<1$. For each $r\ge1$, choose the projection parameters as specified in step (ii) above. Then the DPP approximation satisfies, for all sufficiently small $h>0$,
\begin{flalign*}
 \E\left[\max_{0\le n\le N}|X(t_n)-X^n|^r\right]\le Ch^r. 
\end{flalign*}
Hence the DPPEM method converges strongly with order one for the CEV model.
\end{proposition}

\subsection{AS model}

Consider
\begin{flalign}\label{eq:as} 
\mathrm{d}X_t=\left(\alpha_{-1}X_t^{-1}-\alpha_0+\alpha_1X_t-\alpha_2X_t^{r_A}\right)\mathrm{d}t
+\alpha_3X_t^\rho\,\mathrm{d}B_t, 
\end{flalign}
where all $\alpha_{-1},\alpha_0,\alpha_1,\alpha_2>0$, and $\rho>1$, $r_A>1$, $r_A+1>2\rho$. 

(\romannumeral1) Let $c_0=-\frac{1}{\rho-1}$ and define $ \Phi:(0,\infty)\to(-\infty,0)$ by $\Phi(x)=c_0x^{1-\rho}$ and $\Phi^{-1}(w)=\bigl(\frac{w}{c_0}\bigr)^{-\frac{1}{\rho-1}}$. Then the transformed equation is
\begin{flalign*}
\mathrm{d}W_t=\mu(W_t)\,\mathrm{d}t+\alpha_3\,\mathrm{d}B_t,\quad
W_0=\Phi(X_0),
\end{flalign*}
where
\begin{flalign*}
\mu(w)
=&\,
\alpha_{-1}
\left(\frac{w}{c_0}\right)^{\frac{\rho+1}{\rho-1}}
-\alpha_0
\left(\frac{w}{c_0}\right)^{\frac{\rho}{\rho-1}}
+\alpha_1\left(\frac{w}{c_0}\right)-\alpha_2
\left(\frac{w}{c_0}\right)^{-\frac{r_A-\rho}{\rho-1}}
-\frac{\rho\alpha_3^2}{2}
\left(\frac{w}{c_0}\right)^{-1}.
\end{flalign*}
As $w\uparrow0$, $(b-w)\mu(w)\to -\infty$ because $r_A+1>2\rho$. Hence Assumption \ref{as2} can be satisfied at the finite right
endpoint $b=0$ with arbitrarily large $p_b^\dagger$. Moreover, differentiating $\mu$ shows that $\mu'(w)\to-\infty$ as $w\to-\infty$ and as $w\uparrow0$. Together with continuity, these limits imply that $\mu'$ is bounded above on $(-\infty,0)$. Hence, Assumption~\ref{as1} follows from the mean-value theorem. Therefore Lemma \ref{Lm2.1} yields
\begin{flalign*}
\mathbb P(W_t<0,\ 0\le t\le T)=1,
\end{flalign*}
and $X_t=\Phi^{-1}(W_t)>0$ almost surely.

(\romannumeral2) One has $\gamma_{S,b}=\frac{r_A-1}{\rho-1}$, $\gamma_{I,b}=\frac{2r_A-\rho-1}{\rho-1}$, $\kappa_S=\frac2{\rho-1}$, and $\kappa_I=\frac{\rho+3}{\rho-1}$. Fix $r\ge1$ and set $m:=\frac{\rho}{\rho-1}$, $R:=2r$, $S:=2mr$. Choose $Q>\max\{R,2\}$. Next choose $l_2>0$ sufficiently small such that $l_2\gamma_{S,b}<1$ and $l_2\frac{S}{R}<\frac12$. 
Since $p_b^\dagger$ can be chosen arbitrarily large, choose
$p_b^\dagger$ sufficiently large and then
$p_b^*\in(0,p_b^\dagger)$ such that $Q\gamma_{I,b}<p_b^\dagger$, $S<p_b^\dagger$, $l_2p_b^*>2R$, $l_2\frac{S}{p_b^*}<\frac12$.
It follows that $l_2S (\frac1{p_b^*}+\frac1R)<1$, which is the condition in Lemma \ref{num-moment-inverse} for the numerical inverse moment of order $S$. Choose $c\in(R,l_2p_b^*)$ and then take $P$ sufficiently large so that $c(1-\frac RP)\ge R$ and $\frac cP\,\kappa_S<1$. Set $m_a:=\frac cP$. Then $\delta=\max\{l_2\gamma_{S,b},m_a\kappa_S\}<1$, $\Gamma_X=c$, and 
\begin{flalign*}
\Gamma_X\Big(1-\frac{R}{P}\Big)\ge R \quad\text{and}\quad
\Gamma_X-m_aR\ge R. 
\end{flalign*}

(\romannumeral3) Since $\Phi^{-1}(w)=\left(\frac{w}{c_0}\right)^{-\frac1{\rho-1}}$, its derivative is of order $|w|^{-m}$ near $b=0$, where
$m=\frac{\rho}{\rho-1}$. Hence $|\Phi^{-1}(w)-\Phi^{-1}(v)| \le C|w-v|\bigl(|w|^{-m}+|v|^{-m}\bigr)$ for $ w,v<0$. Theorem \ref{main} gives
\begin{flalign*}
\E\left[
\max_{0\le n\le N}|W(t_n)-W_n|^R
\right]
\le Ch^R.
\end{flalign*}
Moreover, since $S<p_b^\dagger$, Lemma \ref{Lm2.3} gives
\begin{flalign*}
\E\left[\sup_{0\le t\le T}|W_t|^{-S}\right]\le C,
\end{flalign*}
and the condition verified in (\romannumeral2), together with
Lemma \ref{num-moment-inverse}, yield 
\begin{flalign*}
\E\left[
\max_{0\le n\le N}|W_n|^{-S}
\right]\le C.
\end{flalign*}
Therefore, by H\"older's inequality,
\begin{align*}
&\E\left[\max_{0\le n\le N}|X(t_n)-X^n|^r
\right]\\ \le &C\left(
\E\left[\max_{0\le n\le N}|W(t_n)-W_n|^{2r}\right]
\right)^{1/2}\Bigg[\left(\E\Big[\sup_{0\le t\le T}|W_t|^{-S}\Big]\right)^{1/2}
+\left(\E\left[\max_{0\le n\le N}|W_n|^{-S}\right]\right)^{1/2}\Bigg]\\
\le &Ch^r.
\end{align*}

\begin{proposition}
Let all $\alpha_j>0$, $\rho>1$, $r_A>1$, and $r_A+1>2\rho$.  For each $r\ge1$, choose the projection parameters as specified in step (ii) above. Then the DPP approximation satisfies, for all sufficiently small $h>0$,
\begin{flalign*}
\E\left[\max_{0\le n\le N}|X(t_n)-X^n|^r\right]\le Ch^r. 
\end{flalign*}
Hence the DPPEM method converges strongly with order one for the A\"it--Sahalia model.
\end{proposition}

\subsection{Bounded interest-rate model \cite{HHX2027}}

Consider
\begin{flalign*}
\mathrm{d}r_t=(\alpha_0-\beta_0r_t)\,\mathrm{d}t
+\sigma\sqrt{(r_t-r_m)(r_M-r_t)}\,\mathrm{d}B_t,
\quad r_0\in(r_m,r_M). \label{eq:bir} 
\end{flalign*}
This model reduces to the Wright--Fisher model when $r_M=1$ and $r_m=0$. 

(\romannumeral1) Define $\Phi:(r_m,r_M)\to(0,\pi)$ by $\Phi(x)=2\arcsin\sqrt{\frac{x-r_m}{r_M-r_m}}$ with $\Phi^{-1}(w)=r_m+(r_M-r_m)\sin^2(w/2)$. Hence, with $W_t:=\Phi(r_t)$, It\^o's formula yields
\begin{flalign*}
 \mathrm{d}W_t = \mu(W_t)\,\mathrm{d}t + \sigma\,\mathrm{d}B_t,
\quad
W_0=\Phi(r_0)\in(0,\pi),
\end{flalign*}
where
\begin{flalign*}
\mu(w) = \Big(\frac{\alpha_0-\beta_0r_m}{r_M-r_m}-\frac{\sigma^2}{4}\Big) \cot\Big(\frac{w}{2}\Big) -\Big(\frac{\beta_0r_M-\alpha_0}{r_M-r_m}-\frac{\sigma^2}{4}\Big) \tan\Big(\frac{w}{2}\Big), \quad w\in(0,\pi).
\end{flalign*}
At the two finite endpoints,
\begin{flalign*}
\lim_{w\downarrow0}w\mu(w)=2\Big(\frac{\alpha_0-\beta_0r_m}{r_M-r_m}
-\frac{\sigma^2}{4}\Big) \quad\text{and} \quad
\lim_{w\uparrow\pi}(\pi-w)\mu(w)=-2\Big(\frac{\beta_0r_M-\alpha_0}{r_M-r_m}-\frac{\sigma^2}{4}\Big).
\end{flalign*}
Hence, if $\frac{ (\alpha_0-\beta_0r_m) \wedge (\beta_0r_M-\alpha_0) } {\sigma^2(r_M-r_m)} > \frac12$, then Assumption \ref{as2} can be satisfied for every $0<q_a^\dagger<\frac{4(\alpha_0-\beta_0r_m)}{(r_M-r_m)\sigma^2}-2$ and $0<p_b^\dagger<\frac{4(\beta_0r_M-\alpha_0)}{(r_M-r_m)\sigma^2}-2$. Moreover, $\mu'(w)\leq  0$ for $w\in(0,\pi)$, so Assumption \ref{as1} holds. Hence Lemma \ref{Lm2.1} gives
\begin{flalign*}
 \mathbb P(W_t\in(0,\pi),\ 0\le t\le T)=1, 
\end{flalign*}
and therefore $r_t\in(r_m,r_M)$ almost surely.

(\romannumeral2) The exponents are $\gamma_{I,a}=\gamma_{I,b}=3$ and $\gamma_{S,a}=\gamma_{S,b}=2$. Let $\kappa :=4\frac{(\alpha_0-\beta_0r_m) \wedge (\beta_0 r_M-\alpha_0)}{\sigma^2(r_M-r_m)}-2$. Fix $r > 1$ and apply Theorem \ref{main} at transformed moment order $R=r$. Assume $\kappa > 3\max\{r,2\}$. Then one can choose $\max\{r,2\} < Q < \frac{\kappa}{3}$, $3Q < q_a^\dagger < \kappa$, and $3Q < p_b^\dagger < \kappa$. Therefore $\max\{Q\gamma_{I,a},Q\gamma_{S,a}\} = 3Q < q_a^\dagger$ and $\max\{Q\gamma_{I,b},Q\gamma_{S,b}\} = 3Q < p_b^\dagger$.

Since $q_a^\dagger>3Q>2r$ and $p_b^\dagger>3Q>2r$, we may choose $2r<q_a^*<q_a^\dagger$ and $2r<p_b^*<p_b^\dagger$. Then choose the projection exponents $l_1,l_2$ such that $\frac{r}{q_a^*} < l_1 < \frac12$ and $\frac{r}{p_b^*} < l_2 < \frac12$. For $D_h = [h^{l_1},\,\pi-h^{l_2}]$, we have $\delta = \max\{l_1\gamma_{S,a},\, l_2\gamma_{S,b}\} = 2(l_1 \vee l_2) < 1$. Moreover, $l_1q_a^* > r$, $l_2p_b^* > r$, and therefore $\Gamma_X = \min\{l_1q_a^*,\, l_2p_b^*\} > r$ and $\Gamma_X \ge r$. Thus all conditions of Theorem \ref{main} are satisfied at moment order $R=r$.

(\romannumeral3) $\Phi^{-1}(w)=r_m+(r_M-r_m)\sin^2(\frac{w}{2})$ is globally Lipschitz. Consequently, if $\frac{(\alpha_0-\beta_0r_m) \wedge (\beta_0 r_M-\alpha_0)}{\sigma^2(r_M-r_m)}>2$, then every $ 1\le r<\frac13(4\frac{(\alpha_0-\beta_0r_m) \wedge (\beta_0 r_M-\alpha_0)}{\sigma^2(r_M-r_m)}-2)$
is admissible.

\begin{proposition}
Suppose $\frac{(\alpha_0-\beta_0r_m) \wedge (\beta_0 r_M-\alpha_0)}{\sigma^2(r_M-r_m)}>2$ and let $\kappa$ be as defined in step (ii). For each $1\le r<\frac\kappa3$, choose the projection parameters as specified in step (ii). Then the DPPEM approximation satisfies, for all sufficiently small $h>0$,
\begin{flalign*}
\E\left[\max_{0\le n\le N}|r(t_n)-r^n|^r\right]\le Ch^r. 
\end{flalign*}
Hence the scheme is strongly convergent with order one for the bounded interest-rate model.
\end{proposition}

	\section{Numerical examples}
    In this section, we present numerical experiments for the models considered in Section~\ref{sec4} to verify the strong order one convergence of the proposed DPPEM scheme and compare its computational efficiency with the Lamperti backward Euler--Maruyama (LBEM) scheme \cite{NS2014}.  

For comparison, we also consider the LBEM scheme,
\begin{equation*}
    Y_{n+1}^{\rm LBEM}
    =
    Y_n^{\rm LBEM}
    +h\mu(Y_{n+1}^{\rm LBEM})
    +\sigma\Delta B_n,
    \qquad
    Y_0^{\rm LBEM}=\Phi(X_0),
\end{equation*}
with $X_{n+1}^{\rm LBEM} = \Phi^{-1}(Y_{n+1}^{\rm LBEM})$. Throughout the experiments, we take $T=1$ and use the step sizes $2^{-8},\,2^{-9},\,2^{-10},\,2^{-11},\,2^{-12}\,$. The strong error is measured in the discrete pathwise $L^2$ norm
\begin{equation*}
\mathrm{Err}(h) = \left( \frac1M \sum_{m=1}^{M}
\max_{0\le n\le N} \left| X_{\mathrm{ref}}^{(m)}(t_n)
- X_{h}^{(m)}(t_n) \right|^2 \right)^{1/2},
\end{equation*}
where $M=10^4$ Monte Carlo sample paths. Since exact solutions are unavailable for some models, a numerical approximation computed on a finer grid $h_{\mathrm{ref}}=2^{-14}$ is used as the reference solution. The same Brownian paths are used for all schemes and step sizes. The order is the slope of \(\log(\mathrm{Err}(h))\) against \(\log h\) over the five step sizes.

The LBEM update is well defined under the assumptions in \cite{NS2014}. For the CIR and Heston \(3/2\) models, the implicit step can be solved directly, whereas for the other models a nonlinear equation is solved numerically at each time step. In contrast, each DPPEM step requires only an explicit drift evaluation followed by a scalar projection onto \(\mathbb D_h\). The corresponding strong errors, convergence orders, and CPU times are reported below.

\begin{example}{\textbf{CIR model}}

    For this model, we set the projection interval 
\begin{equation*}
\mathbb D_h=[h^{0.45},\,h^{-0.25}]
\end{equation*}
and take the parameters $\kappa=2$, $\theta=0.5$, $\beta=0.4$, $X_0=0.5$. 
For these parameters, $\frac{\kappa\theta}{\beta^2}=6.25>2$ and $\frac{1}{6}\Big( \frac{4\kappa\theta}{\beta^2}-2 \Big) =
\frac{23}{6}>2$, so Proposition~\ref{pro41} applies with $r=2$. The convergence orders of our scheme and LBEM scheme are shown in Figure~\ref{fig:cir_conv}. The numerical results align with our theoretical findings.

\begin{figure}[htbp]
    \centering
    \includegraphics[width=0.8\textwidth]{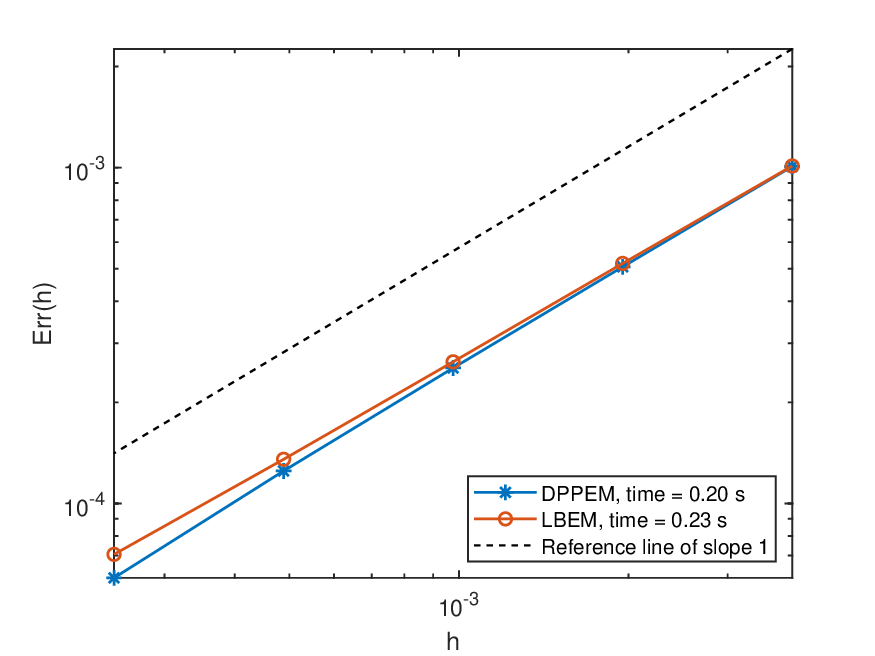}
    \caption{CIR model}
    \label{fig:cir_conv}
\end{figure}

\begin{table}[htbp]
\centering
\caption{Strong errors and convergence orders for the CIR model.}
\label{tab:CIR}
\begin{tabular}{ccccc}
\toprule
$h$ & DPPEM error & DPPEM order & LBEM error & LBEM order \\
\midrule
$3.9062\times10^{-3}$ & $1.00831\times10^{-3}$ & --     & $1.01168\times10^{-3}$ & --     \\
$1.9531\times10^{-3}$ & $5.06525\times10^{-4}$ & 0.9932 & $5.18432\times10^{-4}$ & 0.9645 \\
$9.7656\times10^{-4} $ & $2.53164\times10^{-4}$ & 1.0006 & $2.64089\times10^{-4}$ & 0.9731 \\
$ 4.8828\times10^{-4}$ & $1.24945\times10^{-4}$ & 1.0188 & $1.35336\times10^{-4}$ & 0.9645 \\
$2.4414\times10^{-4}$ & $6.00224\times10^{-5}$ & 1.0577 & $7.05206\times10^{-5}$ & 0.9404 \\
\bottomrule
\end{tabular}
\end{table}
\end{example}

\begin{example}{\textbf{Heston $3/2$ model}}

We choose the model parameters $a_1=2$, $a_2=0.5$, $a_3=0.6$, $X_0=1$, and take the projection interval as
\begin{equation*}
\mathbb D_h=[-h^{-0.25},-h^{0.30}].
\end{equation*}
This choice is compatible with the boundary and negative moment  restrictions required in the theoretical analysis. For these parameters, $\frac{1}{7}\Big( \frac{4a_1}{a_3^2}+2 \Big) = \frac{218}{63}>2$, so Proposition~\ref{pro42} applies with $r=2$. Figure~\ref{heston} displays the corresponding strong orders of the DPPEM scheme and LBEM scheme.

\begin{figure}[htbp]
    \centering
    \includegraphics[width=0.8\textwidth]{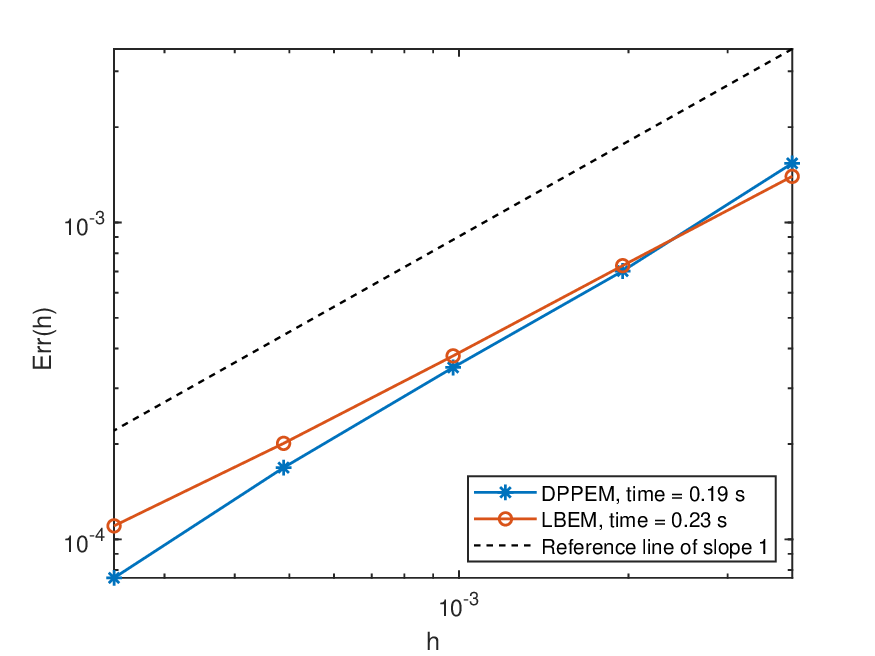}
    \caption{Heston 3/2 model}
    \label{heston}
\end{figure}

\begin{table}[htbp]
\centering
\caption{Strong errors and convergence orders for the Heston $3/2$ model.}
\label{tab:Heston32}
\begin{tabular}{ccccc}
\toprule
$h$ & DPPEM error & DPPEM order & LBEM error & LBEM order \\
\midrule
$3.9062\times10^{-3}$ & $1.53723\times10^{-3}$ & --     & $1.39829\times10^{-3}$ & --     \\
$1.9531\times10^{-3}$ & $7.02693\times10^{-4}$ & 1.1294 & $7.30735\times10^{-4}$ & 0.9362 \\
$9.7656\times10^{-4}$ & $3.48685\times10^{-4}$ & 1.0110 & $3.79337\times10^{-4}$ & 0.9459 \\
$4.8828\times10^{-4}$ & $1.68415\times10^{-4}$ & 1.0499 & $2.00886\times10^{-4}$ & 0.9171 \\
$2.4414\times10^{-4}$ & $7.55801\times10^{-5}$ & 1.1559 & $1.10279\times10^{-4}$ & 0.8652 \\
\bottomrule
\end{tabular}
\end{table}
\end{example}

\begin{example}{\textbf{CEV model}}

The numerical projection interval is chosen as
\begin{equation*}
\mathbb D_h = [h^{0.20},\,h^{-0.40}].
\end{equation*}
The parameters are chosen as $\kappa=2$, $\theta=0.5$, $\beta=0.4$, $\alpha=0.75$, and $X_0=0.1$. Figure~\ref{CEV} reports the observed strong convergence orders of the DPPEM and LBEM schemes. The numerical results support approximately first-order convergence.

\begin{figure}[htbp]
    \centering
    \includegraphics[width=0.8\textwidth]{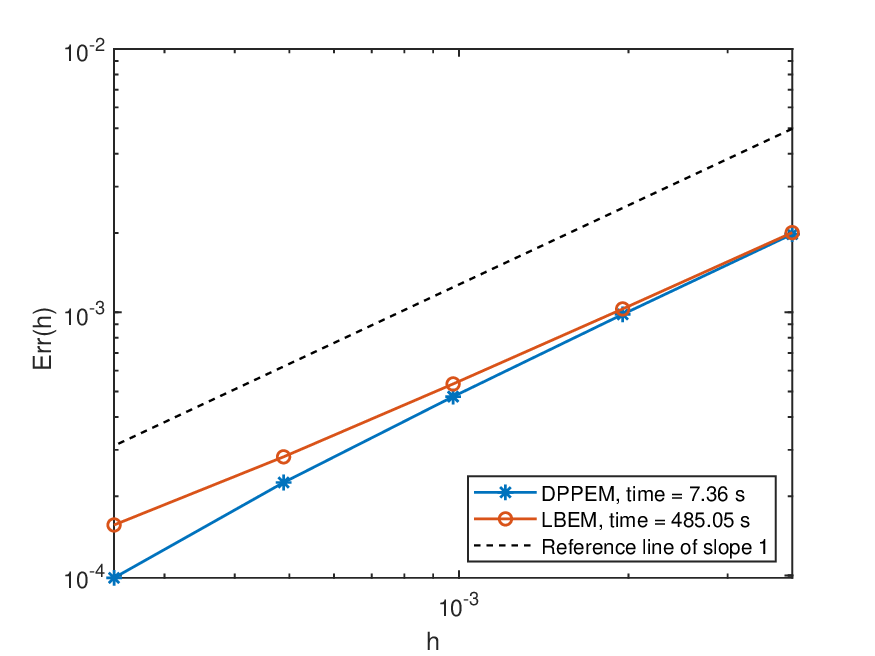}
    \caption{CEV model}
    \label{CEV}
\end{figure}

\begin{table}[htbp]
\centering
\caption{Strong errors and convergence orders for the CEV model.}
\label{tab:CEV}
\begin{tabular}{ccccc}
\toprule
$h$ & DPPEM error & DPPEM order & LBEM error & LBEM order \\
\midrule
$3.9062\times10^{-3}$ & $1.97944\times10^{-3}$ & --     & $2.00585\times10^{-3}$ & --     \\
$1.9531\times10^{-3}$ & $9.81318\times10^{-4}$ & 1.0123 & $1.02918\times10^{-3}$ & 0.9627 \\
$9.7656\times10^{-4}$ & $4.77862\times10^{-4}$ & 1.0381 & $5.33953\times10^{-4}$ & 0.9467 \\
$4.8828\times10^{-4}$ & $2.25480\times10^{-4}$ & 1.0836 & $2.82250\times10^{-4}$ & 0.9197 \\
$2.4414\times10^{-4}$ & $9.79019\times10^{-5}$ & 1.2036 & $1.55509\times10^{-4}$ & 0.8600 \\
\bottomrule
\end{tabular}
\end{table}
\end{example}

\begin{example}{\textbf{AS model}}

We choose the model parameters as $\alpha_{-1}=0.5$, $\alpha_0=0.3$, $\alpha_1=0.2$, $\alpha_2=0.4$, $\alpha_3=0.25$, $\rho=1.5$, $r_A=2.2$, $X_0=1$ and take the projection interval to be
\begin{equation*}
\mathbb D_h = [-h^{-0.24},\,-h^{0.12}].
\end{equation*}
The projection exponents are chosen conservatively so that the stability requirements and the numerical negative moment  conditions used in the theoretical analysis are satisfied. Figure~\ref{AS} shows the resulting strong convergence order, which is consistent with the theoretical orders of two schemes.

\begin{figure}[htbp]
    \centering
    \includegraphics[width=0.8\textwidth]{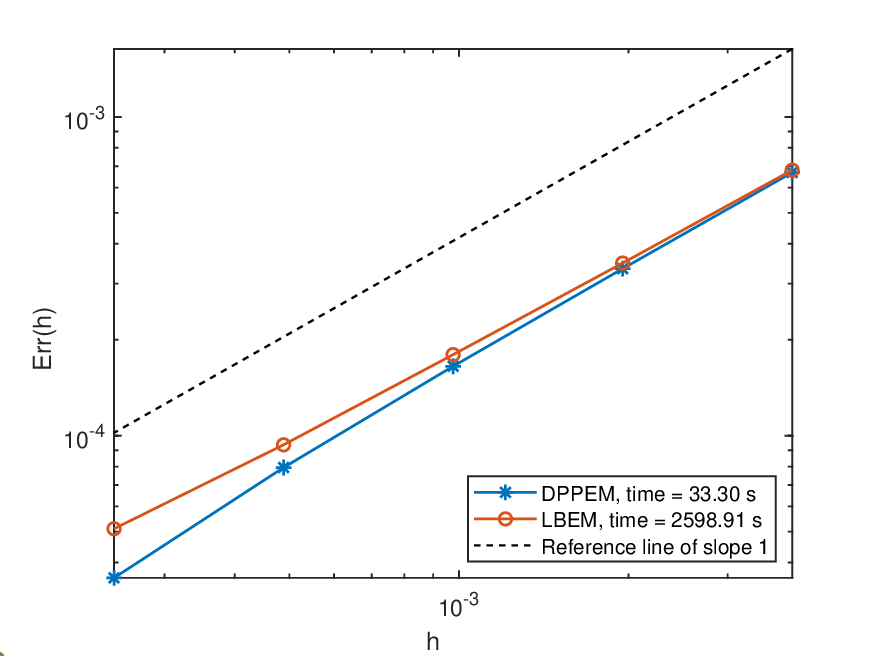}
    \caption{AS model}
    \label{AS}
\end{figure}

\begin{table}[htbp]
\centering
\caption{Strong errors and convergence orders for the A\"it--Sahalia model.}
\label{tab:AS}
\begin{tabular}{ccccc}
\toprule
$h$ & DPPEM error & DPPEM order & LBEM error & LBEM order \\
\midrule
$3.9062\times10^{-3}$ & $6.69650\times10^{-4}$ & --     & $6.80845\times10^{-4}$ & --     \\
$1.9531\times10^{-3}$ & $3.34205\times10^{-4}$ & 1.0027 & $3.47902\times10^{-4}$ & 0.9686 \\
$9.7656\times10^{-4}$ & $1.65116\times10^{-4}$ & 1.0173 & $1.79591\times10^{-4}$ & 0.9540 \\
$4.8828\times10^{-4}$ & $7.94582\times10^{-5}$ & 1.0552 & $9.36272\times10^{-5}$ & 0.9397 \\
$2.4414\times10^{-4}$ & $3.58060\times10^{-5}$ & 1.1500 & $5.10593\times10^{-5}$ & 0.8748 \\
\bottomrule
\end{tabular}
\end{table}
\end{example}

\begin{example}{\textbf{Bounded interest-rate model}}

Consider this model with $\alpha_0=0.15$, $\beta_0=1.5$, $\sigma=0.2$, $r_m=0.02$, $r_M=0.20$, and $r_0=0.08$. Choose
\begin{equation*}
\mathbb D_h = [h^{0.45},\,\pi-h^{0.45}].
\end{equation*}
For the above parameters, $\frac{ (\alpha_0-\beta_0r_m) \wedge (\beta_0r_M-\alpha_0) } {\sigma^2(r_M-r_m)} = \frac{0.12}{0.04\times0.18} >2$,
so the theoretical parameter condition is satisfied. The strong convergence orders are shown in Figure~\ref{BIR}.

\begin{figure}[htbp]
    \centering
    \includegraphics[width=0.8\textwidth]{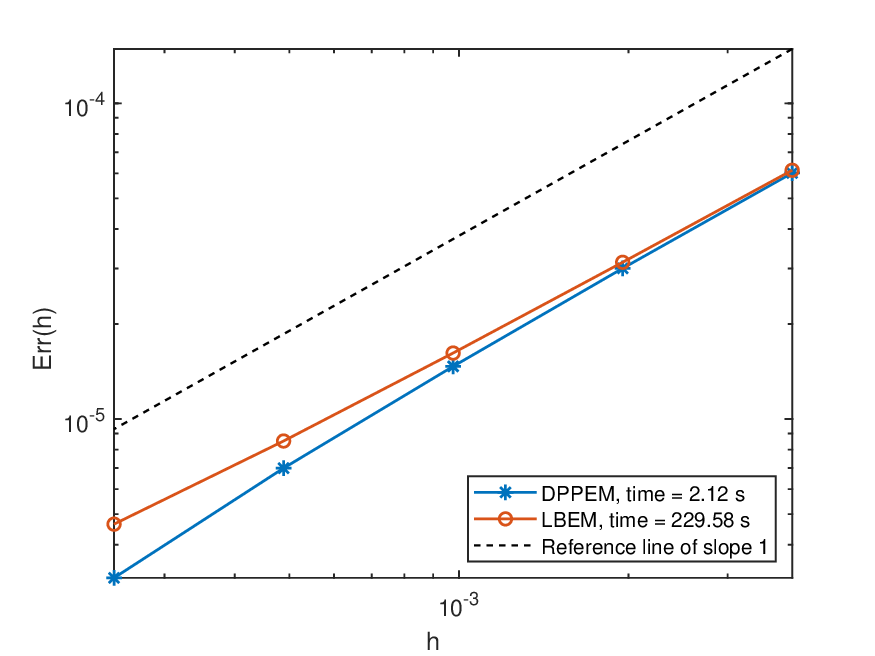}
    \caption{Bounded interest-rate model}
    \label{BIR}
\end{figure}

\begin{table}[htbp]
\centering
\caption{Strong errors and convergence orders for the bounded interest-rate model.}
\label{tab:BIR}
\begin{tabular}{ccccc}
\toprule
$h$ & DPPEM error & DPPEM order & LBEM error & LBEM order \\
\midrule
$3.9062\times10^{-3}$ & $6.01725\times10^{-5}$ & --     & $6.13810\times10^{-5}$ & --     \\
$1.9531\times10^{-3}$ & $3.00206\times10^{-5}$ & 1.0032 & $3.13759\times10^{-5}$ & 0.9681 \\
$9.7656\times10^{-4}$ & $1.46831\times10^{-5}$ & 1.0318 & $1.61859\times10^{-5}$ & 0.9549 \\
$4.8828\times10^{-4}$ & $6.99036\times10^{-6}$ & 1.0707 & $8.52047\times10^{-6}$ & 0.9257 \\
$2.4414\times10^{-4}$ & $3.13991\times10^{-6}$ & 1.1546 & $4.64625\times10^{-6}$ & 0.8749 \\
\bottomrule
\end{tabular}
\end{table}
\end{example}

\begin{table}[htbp]
\centering
\caption{Strong convergence orders and computational costs
of the DPPEM and LBEM schemes.}
\label{tab:order-time}
\begin{tabular}{lccccc}
\toprule
Model
& \multicolumn{2}{c}{Order}
& \multicolumn{2}{c}{CPU time (s)}
& CPU time ratio \\
\cmidrule(lr){2-3}
\cmidrule(lr){4-5}
& DPPEM & LBEM
& DPPEM & LBEM
& LBEM/DPPEM \\
\midrule
CIR
& 1.0160 & 0.9623
& 0.20 & 0.23
& 1.15 \\

Heston $3/2$
& 1.0753 & 0.9192
& 0.19 & 0.23
& 1.21 \\

CEV
& 1.0797 & 0.9245
& 7.36 & 485.05
& 65.90 \\

A\"it--Sahalia
& 1.0523 & 0.9368
& 33.30 & 2598.91
& 78.05 \\

Bounded interest-rate
& 1.0623 & 0.9328
& 2.12 & 229.58
& 108.29 \\
\bottomrule
\end{tabular}
\end{table}

Table~\ref{tab:order-time} summarizes the convergence orders and cumulative solver CPU times. Both schemes exhibit approximately first-order strong convergence over the tested step size. For the CIR and Heston \(3/2\) models, their computational costs are comparable because the LBEM update can be solved directly. For the remaining models, the nonlinear solver  required by LBEM leads to substantially higher computational costs than the fully explicit DPPEM scheme in the present implementation. 
    
	\section{Conclusion}
	In this paper, we proposed a fully explicit DPPEM scheme for scalar SDEs evolving in an interval. Following a Lamperti-type transformation, the scheme combines an explicit EM step with projection onto a step-size-dependent interior interval, preserving the domain pathwise at every time step. Under the stated coefficient assumptions and projection parameter restrictions, we established strong order one convergence for the transformed equation in the $L^r$ norm of the maximum error over the time grid for admissible $r>1$. The analysis accommodated drift singularities near finite endpoints and polynomial growth at infinity. For representative models, the convergence estimates were transferred to the original equations under model specific parameter restrictions. Numerical experiments supported the theoretical convergence rates and showed lower computational costs than the LBEM method for the tested models.
	
	\section*{CRediT authorship contribution statement}
	~~\textbf{Xingwei Hu:} Investigation, Software, Visualization, Writing -- original draft, Writing -- review \& editing. \textbf{Yaozhong Hu:} Funding acquisition, Supervision, Investigation, Validation, Writing -- review \& editing. \textbf{Aiguo Xiao:} Funding acquisition, Supervision, Investigation, Validation, Project administration, Writing -- review \& editing.
	
	\section*{Data Availability}
	~~No external data were used for the research described in the article. 
	
	\section*{Declarations of competing Interests}
	~~The authors declare no competing interests.
	
	\section*{Acknowledgments}
	~~A. Xiao is supported by National Natural Science Foundation of China (No.12471391). X. Hu is supported by the Postgraduate Scientific Research Innovation Project of Hunan Province (No. CX20250930), Xiangtan University (No. XDCX2025Y185), and the  China Scholarship Council under 202508430056. Y. Hu is supported by NSERC Discovery grant RGPIN 2024-05941 and a centennial fund of University of Alberta. 
	
  \bibliographystyle{elsarticle-num}
	\bibliography{refs}

@article{SNS2012,
  author  = {Dereich, Steffen and Neuenkirch, Andreas and Szpruch, Lukasz},
  title   = {An {Euler}-type method for the strong approximation of the {Cox}--{Ingersoll}--{Ross} process},
  journal = {Proc. R. Soc. A},
  volume  = {468},
  number  = {2140},
  pages   = {1105--1115},
  year    = {2012},
  doi     = {10.1098/rspa.2011.0505}
}

@article{AA2005,
	title = {On the discretization schemes for the {CIR} (and {Bessel} squared) processes},
	author = {Alfonsi, Aurélien},
	journal = {Monte Carlo Methods Appl.},
	volume = {11},
	number = {4},
	pages = {355--384},
	year = {2005},
	doi = {10.1515/156939605777438569}
}

@article{AA2013,
	title = {Strong order one convergence of a drift implicit {Euler} scheme: Application to the {CIR} process},
	author = {Alfonsi, Aurélien},
	journal = {Stat. Probab. Lett.},
	volume = {83},
	number = {2},
	pages = {602--607},
	year = {2013},
	doi = {10.1016/j.spl.2012.10.034}
}

@article{AC2020,
  author  = {Cozma, Andrei and Reisinger, Christoph},
  title   = {Strong order 1/2 convergence of full truncation {Euler} approximations to the {Cox}--{Ingersoll}--{Ross} process},
  journal = {IMA J. Numer. Anal.},
  volume  = {40},
  number  = {1},
  pages   = {358--376},
  year    = {2020},
  doi     = {10.1093/imanum/dry067}
}

@article{LCW2025,
  author  = {Liu, Ruishu and Cao, Yulin and Wang, Xiaojie},
  title   = {Unconditionally positivity-preserving explicit {Euler}-type schemes for a generalized {A{\"i}t-Sahalia} model},
  journal = {Numer. Algorithms},
  volume  = {98},
  number  = {2},
  pages   = {651--675},
  year    = {2025},
  doi     = {10.1007/s11075-024-01810-2}
}

@article{JLW2025,
  author  = {Jiang, Yingsong and Liu, Ruishu and Wang, Xiaojie and Zhuo, Jinghua},
  title   = {Unconditionally positivity-preserving approximations of the {A{\"i}t-Sahalia} type model: Explicit {Milstein}-type schemes},
  journal = {Numer. Algorithms},
  volume  = {98},
  number  = {4},
  pages   = {1985--2003},
  year    = {2025},
  doi     = {10.1007/s11075-024-01861-5}
}

@article{hxw2026,
	title={Weak convergence rate of positivity-preserving truncated {Euler}--{Maruyama} method for multi-dimensional stochastic differential equations with positive solutions},
	author={Hu, Xingwei and Dai, Xinjie and Xiao, Aiguo},
	journal={BIT Numer. Math.},
	volume={66},
	number={2},
	pages={19},
	year={2026},
	doi = {10.1007/s10543-026-01117-9}
}

@article{HXW,
	title = {A positivity preserving {Milstein}-type method for stochastic differential equations with positive solutions},
	author = {Hu, Xingwei and Wang, Mengjie and Dai, Xinjie and Yu, Yanyan and Xiao, Aiguo},
	journal = {J. Comput. Appl. Math.},
	volume = {449},
	pages = {115963},
	year = {2024},
	doi = {10.1016/j.cam.2024.115963}
}

@article{DFWM2023,
  author  = {Deng, Shounian and Fei, Chen and Fei, Weiyin and Mao, Xuerong},
  title   = {Positivity-preserving truncated {Euler}--{Maruyama} method for generalised {A{\"i}t-Sahalia}-type interest model},
  journal = {BIT Numer. Math.},
  volume  = {63},
  number  = {4},
  pages   = {59},
  year    = {2023},
  doi     = {10.1007/s10543-023-01000-x}
}

@article{CGW2021,
	title={First order strong convergence of an explicit scheme for the stochastic {SIS} epidemic model},
	author={Chen, Lin and Gan, Siqing and Wang, Xiaojie},
	journal={J. Comput. Appl. Math.},
	volume={392},
	pages={113482},
	year={2021},
	doi = {10.1016/j.cam.2021.113482}
}

@article{YH2021,
	title = {First order strong convergence of positivity preserving logarithmic {Euler}--{Maruyama} method for the stochastic {SIS} epidemic model},
	author = {Yang, Hongfu and Huang, Jianhua},
	journal = {Appl. Math. Lett.},
	volume = {121},
	pages = {107451},
	year = {2021},
	doi = {10.1016/j.aml.2021.107451}
}

@article{LW2023,
	title = {A higher order positivity preserving scheme for the strong approximations of a stochastic epidemic model},
	author = {Liu, Ruishu and Wang, Xiaojie},
	journal = {Commun. Nonlinear Sci. Numer. Simul.},
	volume = {124},
	pages = {107258},
	year = {2023},
	doi = {10.1016/j.cnsns.2023.107258}
}

@article{LC2023,
	title={A positivity preserving {Lamperti} transformed {Euler}--{Maruyama} method for solving the stochastic {Lotka}--{Volterra} competition model},
	author={Li, Yan and Cao, Wanrong},
	journal = {Commun. Nonlinear Sci. Numer. Simul.},
	volume = {122},
	pages = {107260},
	year={2023},
	doi = {10.1016/j.cnsns.2023.107260}
}

@article{LNX2025,
	title = {A strong order 1.5 boundary preserving discretization scheme for scalar {SDEs} defined in a domain},
	author = {Liu, Ruishu and Neuenkirch, Andreas and Wang, Xiaojie},
	journal = {Math. Comp.},
	volume = {94},
	number = {354},
	pages = {1815--1862},
	year = {2025},
	doi = {10.1090/mcom/4014}
}

@article{LG2023,
	title = {Strong and weak convergence rates of logarithmic transformed truncated {EM} methods for {SDEs} with positive solutions},
	author = {Lei, Ziyi and Gan, Siqing and Chen, Ziheng},
	journal = {J. Comput. Appl. Math.},
	volume = {419},
	pages = {114758},
	year={2023},
	doi = {10.1016/j.cam.2022.114758}
}

@article{WG2025,
	title = {An explicit positivity-preserving scheme for the {Heston} 3/2-model with order-one strong convergence},
	author = {Wu, Xiaojuan and Gan, Siqing},
	journal = {Commun. Nonlinear Sci. Numer. Simul.},
	volume = {140},
	pages = {108372},
	year = {2025},
	doi = {10.1016/j.cnsns.2024.108372}
}

@article{CWF2024,
	title={An advanced numerical scheme for multi-dimensional stochastic {Kolmogorov} equations with superlinear coefficients},
	author={Cai, Yongmei and Mao, Xuerong and Wei, Fengying},
	journal={J. Comput. Appl. Math.},
	volume={437},
	pages={115472},
	year={2024},
	doi = {10.1016/j.cam.2023.115472}
}

@article{TX2024,
  author  = {Tang, Yiyi and Mao, Xuerong},
  title   = {The logarithmic truncated {EM} method with weaker conditions},
  journal = {Appl. Numer. Math.},
  volume  = {198},
  pages   = {258--275},
  year    = {2024},
  doi     = {10.1016/j.apnum.2024.01.009}
}

@article{NS2014,
  author  = {Neuenkirch, Andreas and Szpruch, Lukasz},
  title   = {First order strong approximations of scalar {SDEs} defined in a domain},
  journal = {Numer. Math.},
  volume  = {128},
  number  = {1},
  pages   = {103--136},
  year    = {2014},
  doi     = {10.1007/s00211-014-0606-4}
}

@article{PPlogTM,
	title = {Positivity preserving logarithmic {Euler}--{Maruyama} type scheme for stochastic differential equations},
	author = {Yi, Yulian and Hu, Yaozhong and Zhao, Jingjun},
	journal = {Commun. Nonlinear Sci. Numer. Simul.},
	volume = {101},
	pages = {105895},
	year = {2021},
	doi = {10.1016/j.cnsns.2021.105895}
}

@article {zhou_hu_zhao,
	AUTHOR = {Zhou, Hao and Hu, Yaozhong and Zhao, Jingjun},
	TITLE = {Explicit-implicit methods for stochastic
	susceptible-infected-recovered model},
	JOURNAL = {Numer. Algorithms},
	FJOURNAL = {Numer. Algorithms},
	VOLUME = {100},
	YEAR = {2025},
	NUMBER = {3},
	PAGES = {1199--1238},
	ISSN = {1017-1398,1572-9265},
	MRCLASS = {65C30 (60H10 60H35)},
	MRNUMBER = {4970410},
	DOI = {10.1007/s11075-024-01986-7},
}

@article{GGHM2011,
author = {Gray, A. and Greenhalgh, D. and Hu, L. and Mao, X. and Pan, J.},
title = {A Stochastic Differential Equation SIS Epidemic Model},
journal = {SIAM J. Appl. Math.},
volume = {71}, 
number = {3},
pages = {876--902},
year = {2011},
doi = {10.1137/10081856X}
}

@article{HHX2027,
title = {Two novel domain-preserving schemes for a class of stochastic differential equations defined in a domain},
journal = {J. Comput. Appl. Math.},
volume = {492},
pages = {118090},
year = {2027},
issn = {0377-0427},
doi = {10.1016/j.cam.2026.118090},
author = {Xingwei Hu and Yaozhong Hu and Aiguo Xiao}
}

@article{EL2026,
title = {Preserving invariant domains and strong approximation of stochastic differential equations},
journal = {Appl. Numer. Math.},
volume = {221},
pages = {30--45},
year = {2026},
issn = {0168-9274},
doi = {10.1016/j.apnum.2025.11.003},
author = {Utku Erdoğan and Gabriel Lord}
}

@article{MS2007,
author = {Esteban Moro and Henri Schurz },
title = {Boundary Preserving Semianalytic Numerical Algorithms for Stochastic Differential Equations},
journal = {SIAM J. Sci. Comput.},
volume = {29},
number = {4},
pages = {1525--1549},
year = {2007},
doi = {10.1137/05063725X}
}

@article{Niko2015,
  author  = {Halidias, Nikolaos},
  title   = {Construction of positivity preserving numerical schemes for some multidimensional stochastic differential equations},
  journal = {Discrete Contin. Dyn. Syst. Ser. B},
  volume  = {20},
  number  = {1},
  pages   = {153--160},
  year    = {2015},
  doi     = {10.3934/dcdsb.2015.20.153}
}

@misc{brc26,
      title={Analysis of a first-order explicit positivity preserving scheme for a class of scalar {SDEs}}, 
      author={Charles-Edouard Bréhier and David Cohen},
      year={2026},
      eprint={2608.25689},
      archivePrefix={arXiv},
      primaryClass={math.NA}
}

@article{CPC2022,
title = {A new Lagrange multiplier approach for constructing structure preserving schemes, {I}. {Positivity} preserving},
journal = {Comput. Methods Appl. Mech. Engrg.},
volume = {391},
pages = {114585},
year = {2022},
issn = {0045-7825},
doi = {10.1016/j.cma.2022.114585},
author = {Qing Cheng and Jie Shen}
}

@article{CPC2022.2,
author = {Cheng, Qing and Shen, Jie},
title = {A New Lagrange Multiplier Approach for Constructing Structure Preserving Schemes, {II}. {Bound} Preserving},
journal = {SIAM J. Numer. Anal.},
volume = {60},
number = {3},
pages = {970--998},
year = {2022},
doi = {10.1137/21M144877X}
}

@incollection{JOB2015,
  author    = {Ob{\l}{\'o}j, Jan and Spoida, Peter and Touzi, Nizar},
  title     = {Martingale Inequalities for the Maximum via Pathwise Arguments},
  booktitle = {In Memoriam Marc Yor -- S{\'e}minaire de Probabilit{\'e}s XLVII},
  editor    = {Donati-Martin, Catherine and Lejay, Antoine and Rouault, Alain},
  series    = {Lecture Notes in Mathematics},
  volume    = {2137},
  pages     = {227--247},
  publisher = {Springer International Publishing},
  address   = {Cham},
  year      = {2015},
  doi       = {10.1007/978-3-319-18585-9_11}
}

@article{HR2020,
author = {Steven R. Howard and Aaditya Ramdas and Jon McAuliffe and Jasjeet Sekhon},
title = {{Time-uniform Chernoff bounds via nonnegative supermartingales}},
volume = {17},
journal = {Probability Surveys},
number = {none},
publisher = {Institute of Mathematical Statistics and Bernoulli Society},
pages = {257--317},
year = {2020},
doi = {10.1214/18-PS321}
}

@article{HZS2002,
author = {Higham, Desmond J. and Mao, Xuerong and Stuart, Andrew M.},
title = {Strong Convergence of {Euler-Type} Methods for Nonlinear Stochastic Differential Equations},
journal = {SIAM J. Numer. Anal.},
volume = {40},
number = {3},
pages = {1041--1063},
year = {2002},
doi = {10.1137/S0036142901389530}
}

@article{Mehri2021,
  author  = {Mehri, Sima and Scheutzow, Michael},
  title   = {A stochastic {Gronwall} lemma and well-posedness of path-dependent {SDEs} driven by martingale noise},
  journal = {ALEA Lat. Am. J. Probab. Math. Stat.},
  volume  = {18},
  pages   = {193--209},
  year    = {2021},
  doi     = {10.30757/ALEA.v18-09}
}

@article{1977,
     author = {Lenglart, E.},
     title = {Relation de domination entre deux processus},
     journal = {Annales de l'institut Henri Poincar\'e. Section B. Calcul des probabilit\'es et statistiques},
     pages = {171--179},
     year = {1977},
     publisher = {Gauthier-Villars},
     volume = {13},
     number = {2},
     mrnumber = {471069},
     zbl = {0373.60054},
     language = {fr}
}

@article{KAHL2008,
title = {Structure preserving stochastic integration schemes in interest rate derivative modeling},
journal = {Appli. Numer. Math.},
volume = {58},
number = {3},
pages = {284--295},
year = {2008},
issn = {0168-9274},
doi = {10.1016/j.apnum.2006.11.013},
author = {C. Kahl and M. Günther and T. Rossberg}
}

@article{HS2015,
author = {Halidias, N. and Stamatiou, I. S.},
title = {Approximating Explicitly the Mean-Reverting CEV Process},
journal = {J. Probab. Stat.},
volume = {2015},
number = {1},
pages = {513137},
doi = {10.1155/2015/513137},
year = {2015}
}

@article{CIR1985,
  author  = {Cox, John C. and Ingersoll, Jonathan E. and Ross, Stephen A.},
  title   = {A Theory of the Term Structure of Interest Rates},
  journal = {Econometrica},
  volume  = {53},
  number  = {2},
  pages   = {385--407},
  year    = {1985},
  doi     = {10.2307/1911242}
}

@article{HC2025,
title = {Efficient positivity preserving schemes for stochastic complex systems},
journal = {J. Comput. Appl. Math.},
volume = {462},
pages = {116464},
year = {2025},
issn = {0377-0427},
doi = {https://doi.org/10.1016/j.cam.2024.116464},
author = {Can Huang and Huangxin Chen and Qing Cheng and Lijun Chen}
}

@article{BIK2016,
  author  = {Beyn, Wolf-J{\"u}rgen and Isaak, Elena and Kruse, Raphael},
  title   = {Stochastic {C}-stability and {B}-consistency of explicit and implicit {Euler}-type schemes},
  journal = {J. Sci. Comput.},
  volume  = {67},
  number  = {3},
  pages   = {955--987},
  year    = {2016},
  doi     = {10.1007/s10915-015-0114-4}
}

@article{BIK2017,
  author  = {Beyn, Wolf-J{\"u}rgen and Isaak, Elena and Kruse, Raphael},
  title   = {Stochastic {C}-stability and {B}-consistency of explicit and implicit {Milstein}-type schemes},
  journal = {J. Sci. Comput.},
  volume  = {70},
  number  = {3},
  pages   = {1042--1077},
  year    = {2017},
  doi     = {10.1007/s10915-016-0290-x}
}

@article{THX2026,
  author  = {Tao, Pengfei and Hu, Xingwei and Xiao, Aiguo},
  title   = {Positivity-preserving {Lamperti}-transformed truncated {Euler}--{Maruyama} method for the stochastic delay {SIRs} model},
  journal = {Int. J. Comput. Math.},
  year    = {2026},
  doi     = {10.1080/00207160.2026.2723792}
}
\end{document}